\documentclass[a4paper,12pt,pagesize]{scrartcl}
\pdfoutput=1

\usepackage[utf8]{inputenc}
\usepackage[english,ngerman]{babel}
\usepackage[T1]{fontenc}
\usepackage{lmodern}
\usepackage{microtype}
\usepackage{amsmath}
\usepackage{amssymb}
\usepackage{hyperref}
\usepackage{hypcap}
\usepackage[hyperref,amsmath,thmmarks]{ntheorem}
\usepackage{tikz}

\newcommand{\R}{\mathbb R}
\newcommand{\Q}{\mathbb Q}
\newcommand{\Z}{\mathbb Z}

\newcommand{\F}{\mathbb F}

\newcommand{\bP}{\mathbb P}

\newcommand{\cA}{\mathcal A}

\newcommand{\cC}{\mathcal C}
\newcommand{\cD}{\mathcal D}
\newcommand{\cF}{\mathcal F}

\newcommand{\cP}{\mathcal P}

\newcommand{\cX}{\mathcal X}
\newcommand{\cY}{\mathcal Y}
\newcommand{\cZ}{\mathcal Z}
\newcommand{\caL}{\mathcal L}
\newcommand{\cI}{\mathcal I}

\DeclareMathOperator{\id}{id}
\DeclareMathOperator{\im}{im}

\DeclareMathOperator{\Aut}{Aut}

\DeclareMathOperator{\PGl}{PGl}
\DeclareMathOperator{\Sl}{Sl}
\DeclareMathOperator{\PSl}{PSl}
\DeclareMathOperator{\St}{St}
\DeclareMathOperator{\Sp}{Sp}

\DeclareMathOperator{\lk}{lk}
\DeclareMathOperator{\Lk}{Lk}
\DeclareMathOperator{\st}{st}

\newcommand{\coloneq}{\mathrel{\mathop :}=}

\theoremstyle{plain}
\newtheorem{Definition}{Definition}[section]
\newtheorem{Theorem}[Definition]{Theorem}
\newtheorem{Lemma}[Definition]{Lemma}
\newtheorem{Proposition}[Definition]{Proposition}
\newtheorem{Corollary}[Definition]{Corollary}

\theoremclass{Theorem}
\theoremstyle{break}

\theoremstyle{nonumberplain}
\theorembodyfont{\rm}
\newtheorem{Remark}{Remark}
\newtheorem{Remarks}{Remarks}
\newtheorem{Examples}{Examples}
\theoremsymbol{\ensuremath\Box}
\newtheorem{Proof}{Proof}

\addtokomafont{paragraph}{\rm\bf}
\addtokomafont{subparagraph}{\rm\it}

\title{A geometric construction of panel-regular lattices in buildings of types $\tilde A_2$ and $\tilde C_2$}
\author{Jan Essert}
\date{\today}

\xdefinecolor{myblue}{rgb}{0.1,0.1,0.6}
\xdefinecolor{mygreen}{rgb}{0.1,0.4,0.3}

\hypersetup{
pdftitle={A geometric construction of panel-regular lattices in buildings of types ~A2 and ~C2},
pdfauthor={Jan Essert},
pdfkeywords={Affine buildings, Lattices, Exotic buildings, Group theory, Complexes of groups},
colorlinks=true,
linkcolor=myblue,
citecolor=mygreen,
urlcolor=myblue
}

\begin{document}
\selectlanguage{english}
\maketitle

\begin{abstract}
    Using Singer polygons, we construct locally finite affine buildings of types $\tilde A_2$ and $\tilde C_2$ which admit uniform lattices acting regularly on panels. This construction produces very explicit descriptions of these buildings as well as very short presentations of the lattices. All but one of the $\tilde C_2$-buildings are necessarily exotic. To the knowledge of the author, these are the first presentations of lattices in buildings of type $\tilde C_2$. Integral and rational group homology for the lattices is also calculated.
\end{abstract}

\addtocounter{section}{-1}
\section{Introduction}

A \emph{uniform lattice} is a discrete and cocompact subgroup in a locally compact topological group. In recent years, the investigation of lattices in the automorphism groups of locally finite polyhedral complexes with non-positive curvature has been an active topic of research. Many open problems and an overview of the current situation can be found in the survey \cite{FHT:PAG:09}.

This text describes a new geometric construction method producing buildings of type $\tilde A_2$ and $\tilde C_2$ with a uniform lattice in their full automorphism groups. The advantages of this construction are very simple presentations for the lattices, very explicit descriptions of the buildings and a very good understanding how these lattices act on the building. Using the latter fact, we can calculate the group homology of these lattices.

The lattices act regularly, that is transitively and freely, on panels of the same type of the building. It is known that there are countably many buildings of types $\tilde A_2$ and $\tilde C_2$ coming from algebraic groups over local fields, as well as uncountably many other so-called \emph{exotic} buildings. For the buildings of type $\tilde C_2$, all but possibly one of the buildings we construct are certainly exotic.

In the case of $\tilde A_2$, this is not known yet, except for a single lattice which we can realise explicitly in $\Sl_3(\F_2(\!(t)\!))$. If the building $X$ we construct is classical, then there is an associated algebraic group $G\leq \Aut(X)$. By a classical result of Tits in \cite{Tit:BsT:74}, the algebraic group $G$ is always cocompact in $\Aut(X)$ and it is hence conceivable that our lattice (or a finite-index subgroup) is already contained in $G$. Then Margulis' Arithmeticity \cite{Mar:DSS:91} asserts that the lattice is \emph{arithmetic}, that is, comes from an algebraic construction. Our construction would then provide simple presentations for arithmetic lattices, along with a very explicit description of the action on the building.

\paragraph{Details} For the construction, we use \emph{Singer polygons}, generalised polygons with a point-regular automorphism group called a \emph{Singer group}, to construct small complexes of groups. The local developments of these complexes of groups are cones over the generalised polygons we started with. These cones are automatically non-positively curved in a natural metric, so the complexes of groups are developable by a theorem of Bridson and Haefliger in \cite{BH:NPC:99}. Then, by a recognition theorem by Charney and Lytchak in \cite{CL:MC:01}, we know that their universal covers are buildings. The fundamental groups of these complexes of groups are then uniform lattices on two-dimensional locally finite affine buildings.

In principle, our construction might even extend to lattices in buildings of type $\tilde G_2$, but unfortunately, there are no known Singer hexagons.

\paragraph{Previous constructions} There are several previous geometric constructions of cham\-ber-tran\-si\-tive lattices in buildings which led to the classification of all cham\-ber-tran\-si\-tive lattices in classical buildings by Kantor, Liebler and Tits in \cite{KLT:CTL:87}, among these three lattices in buildings of type $\tilde A_2$.

Cartwright, Mantero, Steger and Zappa construct vertex-transitive lattices on buildings of type $\tilde A_2$ in \cite{CMSZ1, CMSZ2}. They prove that some of the lattices they construct are contained in $\Sl_3(\F_q(\!(t)\!))$, while others correspond to exotic buildings.

In \cite{Ron:TG:84}, Ronan constructs possibly exotic buildings of type $\tilde A_2$ admitting a lattice acting regularly on vertices of the same type. Except for the single lattice mentioned above, it is not clear whether any of our lattices is commensurable to the lattices constructed by Cartwright-Mantero-Steger-Zappa or by Ronan.

\paragraph{} For buildings of type $\tilde C_2$, there is a free construction by Ronan in \cite{Ron:CBR:86} which produces very unstructured examples of exotic buildings, but which does not give any control over the automorphism groups. In \cite{Kan:GPS:86}, Kantor gives a construction of exotic buildings of type $\tilde C_2$ admitting uniform lattices acting freely on vertices of the building. However, there are no presentations of these lattices given. Again, we do not know whether the lattices we construct here are commensurable to the examples given by Kantor.

\subsection{Buildings of type \texorpdfstring{$\tilde A_2$}{\textasciitilde A2}}

For buildings of type $\tilde A_2$, we obtain the following result. Fix three projective planes of order $q$ along with three Singer groups $S_i$, $i\in\{1,2,3\}$, that is, groups acting regularly on points (and hence on lines) of these planes. For each of these planes, fix a point $p_i$ and a line $l_i$ and write
\[
    D_i = \{ d\in S_i: p_i\text{ is incident to } d(l_i) \},
\]
these sets are called \emph{difference sets}. Finally, write $J=\{0,1,\ldots,q\}$ and fix three bijections $d_i:J\rightarrow D_i$.

\begin{Theorem}\label{th:intr_a2_lattices}
    There is a building $X$ such that the group
    \[
    \Gamma_1 = \Bigl\langle S_1,S_2,S_3 \,\Big| \begin{array}{c}\text{ all relations in the groups }S_1,S_2,S_3,\\ d_1(j)d_2(j)d_3(j) = d_1(j')d_2(j')d_3(j') \quad\forall j,j'\in J\end{array}\Bigr\rangle
    \]
is a uniform lattice in the full automorphism group of the building. The set of all chambers adjacent to any panel is a fundamental domain for the action. Moreover, every panel-regular lattice on a building of type $\tilde A_2$ arises in this way and admits a presentation of the above form with $d_1(0)d_2(0)d_3(0)=1$.
\end{Theorem}

\noindent A very simple special case arises as follows:

Classical projective planes admit cyclic Singer groups. Denote the generators of the three Singer groups $S_i$ by $\sigma_i$, and write
\[
\Delta_i=\{ \delta \in \Z / |S_i| : \sigma_i^\delta \in D_i\},
\]
these are difference sets in the classical sense.

\begin{Theorem}\label{th:intr_cyclic_lattices}
    For any prime power $q$, for any three classical difference sets $\Delta_1$, $\Delta_2$, $\Delta_3$ containing 0, and for any bijections $\delta_\alpha:J\rightarrow \Delta_\alpha$ satisfying $\delta_\alpha(0)=0$, the group $\Gamma_2$ with presentation
	\[
		\Gamma_2 = \langle \sigma_1, \sigma_2, \sigma_3 \,|\, \sigma_1^{q^2+q+1}=\sigma_2^{q^2+q+1}=\sigma_3^{q^2+q+1}=1,\, \sigma_1^{\delta_1(j)}\sigma_2^{\delta_2(j)}\sigma_3^{\delta_3(j)} = 1 \quad\forall j \in J\rangle
\]
	is a uniform lattice in a building of type $\tilde A_2$.
\end{Theorem}

\begin{Examples}
    Two simple examples are:
\begin{align*}
    \Lambda &= \langle \sigma_1,\sigma_2,\sigma_3 \,|\, \sigma_1^7=\sigma_2^7=\sigma_3^7 = \sigma_1\sigma_2\sigma_3 = \sigma_1^3\sigma_2^3\sigma_3^3 = 1\rangle\\
    \Lambda' &= \langle \sigma_1,\sigma_2,\sigma_3 \,|\, \sigma_1^{13}=\sigma_2^{13}=\sigma_3^{13} = \sigma_1\sigma_2^3\sigma_3^9 = \sigma_1^3\sigma_2^9\sigma_3 =\sigma_1^9\sigma_2\sigma_3^3 = 1\rangle
\end{align*}
Further examples can easily be constructed using the list of difference sets found in the La Jolla Difference Set Repository \cite{LaJolla}.
\end{Examples}

If we start from cyclic Singer groups, we obtain very explicit descriptions of the associated buildings. Using these, we can give explicit descriptions of the spheres of radius 2 in these buildings, which lead to an interesting incidence structure called a \emph{Hjelmslev plane of level two}. Using a result of Cartwright-Mantero-Steger-Zappa in \cite{CMSZ2}, we show that the lattices from Theorem \ref{th:intr_cyclic_lattices} cannot belong to the building associated to $\Sl_3(\Q_p)$ if $\Delta_1=\Delta_2=\Delta_3$.

For the homology groups, we obtain the following result.

\begin{Theorem} If $\Gamma_2$ is a lattice as in Theorem \ref{th:intr_cyclic_lattices}, then
    \[
        H_j(\Gamma_2;\Z) \cong\begin{cases}
            \Z & j=0\\
            \ker(\cD) & j=1 \\
            \Z^q & j=2 \\
            (\Z/(q^2+q+1))^3 & j\geq 3 \text{ odd} \\
            0 & \text{ else,}
        \end{cases}
    \]
    where $\cD:(\Z/(q^2+q+1))^3\rightarrow (\Z/(q^2+q+1))^q$ is given by the matrix $(\delta_j(i))_{i,j}$.

    In addition
    \[
        H_2(\Gamma_i;\Q)=\Q^q,\qquad H_j(\Gamma_i;\Q)=0 \text{ for $j\neq 0,2$}
    \]
    for any lattice $\Gamma_i$ as in Theorems \ref{th:intr_a2_lattices} or \ref{th:intr_cyclic_lattices}.
\end{Theorem}

\subsection{Buildings of type \texorpdfstring{$\tilde C_2$}{\textasciitilde C2}}

For buildings of type $\tilde C_2$, we give two different constructions of panel-regular lattices, depending on the types of panels the lattice acts regularly on. We use \emph{slanted symplectic quadrangles} as defined in \cite{GJS:SSQ:94} for the vertex links. Since almost all of these quadrangles are exotic, we necessarily construct exotic buildings admitting panel-regular lattices.

Starting from a slanted symplectic quadrangle of order $(q-1,q+1)$, where $q>2$ is a prime power, we write $J=\{0,1,\ldots,q+1\}$ and consider a Singer group $S$ acting regularly on points of the quadrangle. Fix the set of lines $L$ through a point and a bijective enumeration function $\lambda:J\rightarrow L$. All line stabilisers $S_l$ are isomorphic to $\Z/q$, and we fix isomorphisms $\psi_j:\Z/q\rightarrow S_{\lambda(j)}$. We repeat this construction for a second quadrangle of order $(q-1,q+1)$ and obtain a Singer group $S'$ and functions $\lambda'$ and $\{\psi'_j\}_{j\in J}$.

\begin{Theorem}\label{th:intr_c2}
    The finitely presented groups
    \begin{align*}
    \Gamma_1 &= (S * S')/ \langle [S_{\lambda(j)}, S'_{\lambda'(j)}] : j\in J\rangle \\
    \intertext{as well as}
    \Gamma_2 &= (S * S' * \langle c\rangle) / \langle c^{q+2},\, c^j \psi_j(x) c^{-j} \psi_j'(x)^{-1} : j\in J, x\in \Z/q\rangle
    \end{align*}
    are uniform panel-regular lattices in buildings of type $\tilde C_2$.

    \begin{itemize}
        \item If $q$ is an odd prime power, then $S$ and $S'$ are three-dimensional Heisenberg groups over $\F_q$. If $q$ is an odd prime, the above presentations can be made more explicit by writing out presentations for $S$ and $S'$.
        \item If $q$ is even, then $S$ and $S'$ are isomorphic to the additive groups of $\F_q^3$. The set of lines can be identified with $\bP\F_q^2\sqcup \{0\}$, the stabilisers have the following structure:
            \begin{itemize}
                \item For $[a:b]$ in the projective plane $\bP\F_q^2$, the stabiliser $S_{[a:b]}$ is the $\F_q$-subspace of $S=\F_q^3$ spanned by $(a,b,0)^T$. 
                \item The stabiliser $S_0$ is the $\F_q$-subspace spanned by $(0,0,1)^T$.
            \end{itemize}
    \end{itemize}
\end{Theorem}

\noindent If $q>3$, then the associated buildings are necessarily exotic. In any case, the above presentations imply very simple descriptions of these buildings. Finally, we calculate group homology.

\begin{Theorem}
    If $\Gamma_i$ is any of the two lattices constructed in Theorem \ref{th:intr_c2}, we have $H_j(\Gamma_i;\Q)=0$ for $j\neq 0$. In addition, for the first type of lattices we have
    \[
        H_1(\Gamma_1;\Z)\cong (\Z/q)^6,\qquad H_2(\Gamma_1;\Z) \cong H_2(S) \oplus H_2(S').
    \]
\end{Theorem}

\subsection{Overview}

This paper is structured as follows: In section \ref{sec:complexes_of_groups}, we introduce complexes of groups and the concept of developability, the universal cover and local developments, which we need for our constructions. We briefly introduce two-dimensional affine buildings in section \ref{sec:buildings} and state the metric recognition theorem by Charney and Lytchak. Singer polygons and, in particular, the explicit construction of slanted symplectic quadrangles and their Singer groups are treated in section \ref{sec:singer_polygons}. The construction of a complex of groups with cones over Singer polygons as local developments, which is central to this paper, is detailed in section \ref{sec:local_cplx_of_groups_construction}. We construct the panel-regular lattices in buildings of type $\tilde A_2$ and $\tilde C_2$ in sections \ref{sec:a2} and \ref{sec:c2}, respectively. Finally, we calculate integral and rational group homology of these lattices in section \ref{sec:group_homology}.

\paragraph{} The foundations of this paper were laid while the author visited Benson Farb at the University of Chicago. The author wishes to thank Benson Farb for his hospitality and the introduction to the topic of lattices in polyhedral complexes, as well as for asking the questions which led to the investigation of this subject. In addition, we thank Linus Kramer for encouraging the author to pursue this topic further as well as for a lot of additional help. Many people have contributed expert knowledge on various pieces of this subject, we thank Uri Bader, Pierre-Emmanuel Caprace, Theo Grundhöfer, Hendrik van Maldeghem, Markus Stroppel, Koen Thas and Alain Valette for their help.

The author was supported by the \selectlanguage{ngerman}\emph{Graduiertenkolleg: \glqq Analytische Topologie und Me\-ta\-geo\-me\-trie\grqq}\selectlanguage{english} while working on this topic. This work is part of the author's doctoral thesis at the Universität Münster.

\section{Complexes of Groups}\label{sec:complexes_of_groups}

To fix the notation, we will collect some facts about complexes of groups, but our treatment will not be complete. For a reference, see, of course, \cite[Chapter III.$\cC$]{BH:NPC:99}. We will also mostly use the notation of this book, except denoting vertices by $v,w$ instead of $\sigma,\tau$ and changing the notation for the presentations of fundamental groups slightly.

\subsection{Definitions}

\begin{Definition}
	A \emph{small category without loops (scwol)} is a set $\cX$ which is the disjoint union of a vertex set $V(\cX)$ and an edge set $E(\cX)$ together with initial and terminal vertex maps $i,t:E(\cX)\rightarrow V(\cX)$. We denote by $E^{(2)}(\cX)$ the set of pairs $(a,b)\in E(\cX)\times E(\cX)$ such that $i(a)=t(b)$ and require the existence of a composition map $E^{(2)}(\cX)\rightarrow E(\cX)$, $(a,b)\mapsto ab$, satisfying certain additional conditions, which will not be made explicit here.
\end{Definition}

\noindent Associated to any scwol $\cX$ there is the \emph{geometric realisation $|\cX|$}: a piecewise Euclidean complex whose vertices are in bijection to $V(\cX)$ and whose $k$-simplices correspond to $k$-chains of composable edges in $E(\cX)$. For any vertex $v\in V(\cX)$ we denote the \emph{closed star} of $v$ in $|\cX|$ by $\St(v)$. The \emph{open star} or residue is denoted by $\st(v)$.

\begin{Definition}
	Let $\cY$ be a scwol. A \emph{complex of groups} $G(\cY)=(G_v,\psi_a,g_{a,b})$ over $\cY$ is given by the following data:
	\begin{enumerate}
		\item A local group $G_v$ for each vertex $v\in V(\cY)$.
		\item A monomorphism $\psi_a:G_{i(a)}\rightarrow G_{t(a)}$ for each edge $a\in E(\cY)$.
		\item A twisting element $g_{a,b}\in G_{t(a)}$ for each pair of composable edges $(a,b)\in E^{(2)}(\cY)$.
	\end{enumerate}
	with certain compatibility conditions. Throughout this paper we will only consider finite vertex groups.
\end{Definition}

\noindent If a group $G$ acts on a scwol $\cX$, there is a natural way to construct a quotient complex of groups $G\backslash\!\backslash \cX$ over the quotient scwol $G\backslash\cX$. This construction assigns the vertex stabilisers as local groups. A complex of groups is called \emph{developable} if it is isomorphic to a quotient complex arising in this fashion.

The notation in the following definition differs from the one in \cite{BH:NPC:99}. In addition, we simplified the presentation slightly. It is very easy to see that this definition is equivalent to the one in the aforementioned book.

\begin{Definition}\label{def:fundamental_group}
	Let $G(\cY)$ be a complex of groups. We consider the associated graph with vertex set $V(\cY)$ and edge set $E(\cY)$. Choose a maximal tree $T$ in this graph. The \emph{fundamental group $\pi_1(G(\cY),T)$ of $G(\cY)$ with respect to $T$} is generated by the set
	\[
	\coprod_{v\in V(\cY)}G_v \sqcup \{ k_a : a\in E(\cY) \}
	\]
	subject to the relations in the groups $G_v$ and to
	\begin{align*}
		k_a k_b &= g_{a,b} k_{ab}& &\forall (a,b)\in E^{(2)}(\cY)\\
		\psi_a(g) &= k_agk_a^{-1}& &\forall g\in G_{i(a)}\\
		k_a &= 1& &\forall a\in T.
	\end{align*}
\end{Definition}

\noindent It can be shown that fundamental groups associated to different maximal trees are isomorphic. We will hence omit $T$ usually and we will write simply $\pi_1(G(\cY))$. The importance of complexes of groups for the construction of lattices comes from the following assertion.

\begin{Proposition}
    Associated to each developable complex of groups $G(\cY)$, there is a simply connected scwol $\cX$, called the \emph{universal cover} or \emph{development} of $\cY$ on which $\pi_1(G(\cY))$ acts such that $G(\cY)$ is isomorphic to the complex of groups $\pi_1(G(\cY))\backslash\!\backslash \cX$.
\end{Proposition}

\paragraph{Construction} We will need the precise construction of the universal cover $\cX$ later on. Abbreviate $\Gamma\coloneq \pi_1(G(\cY))$. Then write $\varphi_v:G_v\rightarrow \Gamma$ for the canonical homomorphisms, which are injective if and only if $\cY$ is developable. The universal cover $\cX$ is given by
\begin{align*}
	V(\cX)&\coloneq \Bigl\{ \bigl(g\varphi_v(G_v), v\bigr) : g\in \Gamma, v\in V(\cY)\Bigr\} \\
	E(\cX)&\coloneq \Bigl\{ \bigl(g\varphi_{i(a)}(G_{i(a)}),a\bigr): g\in \Gamma, a\in E(\cY)\Bigr\}
\end{align*}
with initial and terminal vertex maps
\begin{align*}
	i\bigl( \bigl(g\varphi_{i(a)}(G_{i(a)}),a\bigr)\bigr) &\coloneq \bigl(g\varphi_{i(a)}(G_{i(a)}),i(a)\bigr)\\
	t\bigl( \bigl(g\varphi_{i(a)}(G_{i(a)}),a\bigr)\bigr) &\coloneq \bigl(gk_a^{-1}\varphi_{t(a)}(G_{t(a)}),t(a)\bigr).
\end{align*}
Composition is defined as follows
\[
(g\varphi_{i(a)}(G_{i(a)}),a)(h\varphi_{i(b)}(G_{i(b)}),b) = (h\varphi_{i(b)}(G_{i(b)}),ab).
\]
These edges are composable if $a$, $b$ are composable and $g\varphi_{i(a)}(G_{i(a)})=hk_b^{-1}\varphi_{i(a)}(G_{i(a)})$.

\begin{Remark}
	If $\cY$ is finite, the fundamental group acts cocompactly on the geometric realisation of $\cX$, all vertex stabilisers are finite. If the scwol $\cY$ is locally finite, so is the universal cover $\cX$. The automorphism group of the geometric realisation $|\cX|$ is then a locally compact group with compact open stabilisers. The fundamental group $\pi_1(G(\cY))$ is discrete and hence a \emph{uniform lattice} in $\Aut(|\cX|)$.
\end{Remark}

\subsection{Developability and non-positive curvature}

\noindent However, not all complexes of groups are developable. If the polyhedral complex $|\cY|$ is endowed with a locally Euclidean metric then there is a criterion for developability assuming non-positive curvature on a local development, as defined in \cite[III.$\cC$.4.20]{BH:NPC:99}.

To ease reading of this paper, we duplicate the construction of the local development here. The reader is assumed to be familiar with the non-positive curvature conditions CAT(0) and CAT(1). A very good introduction can be found in \cite[II]{BH:NPC:99}.

Two auxiliary concepts --- the lower and the upper link --- are required for the definition of the local development.

\begin{Definition}\label{def:upper_link}
	For a vertex $v\in\cY$ in a complex of groups $G(\cY)$ we define the scwol $\Lk_{\tilde v}(\cY)$ as follows:
	\begin{align*}
		V(\Lk_{\tilde v}(\cY)) &= \bigl\{ \bigl(g\psi_a(G_{i(a)}),a\bigr) : a\in E(\cY), t(a)=v, g\psi_a(G_{i(a)})\in G_v/\psi_a(G_{i(a)}) \bigr\} \\
		E(\Lk_{\tilde v}(\cY)) &= \bigl\{ \bigl(g\psi_{ab}(G_{i(b)}),a,b\bigr) : (a,b) \in E^{(2)}(\cY), t(a)=v, g\psi_{ab}(G_{i(b)}) \in G_v/\psi_{ab}(G_{i(b)}) \bigr\}
	\end{align*}
	The maps $i,t : E(\Lk_{\tilde v}(\cY))\rightarrow V(\Lk_{\tilde v}(\cY))$ are defined by
	\begin{align*}
		i\bigl((g\psi_{ab}G_{i(b)},a,b)\bigr) &= \bigl(g\psi_{ab}(G_{i(b)}),ab\bigr),\\
		t\bigl((g\psi_{ab}G_{i(b)},a,b)\bigr) &= \bigl(gg^{-1}_{a,b}\psi_{ab}(G_{i(b)}),a\bigr).
	\end{align*}
	Since all scwols in this paper are at most two-dimensional and hence the scwols $\Lk_{\tilde v}(\cY)$ are at most one-dimensional, we omit the definition of the composition of edges.
\end{Definition}

\noindent We also need the definition of the lower link.

\begin{Definition}
	Let $v\in V(\cY)$ be a vertex in a scwol $\cY$. The \emph{lower link} $\Lk^v(\cY)$ is given by
\begin{align*}
	V(\Lk^v(\cY)) &= \{ a\in E(\cY) : i(a)=v \} \\
	E(\Lk^v(\cY)) &= \{ (a,b) \in E^{(2)}(\cY) : i(b)=v \}
\end{align*}
	with $i(a,b)\coloneq b$, $t(a,b)\coloneq ab$. Again, we do not define composition of edges.
\end{Definition}

\noindent We will not define the join of scwols formally, it is enough to know that the geometric realisation of a join is affinely isometric to the (simplicial) join of the geometric realisations of the two scwols involved.

\begin{Definition}\label{def:local_development}
	For any vertex $v\in V(\cY)$ of a complex of groups $G(\cY)$, we define the \emph{local development of the complex of groups $G(\cY)$ at the vertex $v$} to be the scwol $\cY(\tilde v)\coloneq \Lk^v(\cY) * \{v\} * \Lk_{\tilde v}(\cY)$.

    Notation: We write $\St(\tilde v)=|\cY(\tilde v)|$ for the geometric realisation of $\cY(\tilde v)$ and $\st(\tilde v)$ for the open star of $\tilde v$ in $\St(\tilde v)$.
\end{Definition}

There is a natural projection of $\St(\tilde v)$ onto the closed star $\St(v)$ of the vertex $v$. Assume that the geometric realisation $|\cY|$ of $\cY$ is endowed with a locally Euclidean metric. Then $\St(v)$ is endowed with the induced metric, and we can endow $\St(\tilde v)$ and its subspace $\st(\tilde v)$ with a locally Euclidean metric such that the aforementioned projection is isometric on every simplex.

\begin{Proposition}[III.$\cC$.4.11 in \cite{BH:NPC:99}]\label{prop:local_isometry}
	Assume the complex of groups $G(\cY)$ is developable and let $\cX$ be the universal cover. Let $v\in V(\cY)$ and choose a preimage $\bar v\in\cX$ of $v$. Then there is a $G_v$-equivariant isometry $\st(\tilde v) \rightarrow \st(\bar v)$.
\end{Proposition}

\noindent Hence we have a precise understanding of the local structure of the universal cover, even though the global structure may be complicated.

\begin{Definition}
	If $\lvert\cY\rvert$ can be endowed with a locally Euclidean metric such that for every vertex $v$ of $G(\cY)$ the metric space $\st(\tilde v)$ is CAT(0), we say that the complex of groups is \emph{non-positively curved}.
\end{Definition}

\begin{Theorem}[Bridson-Haefliger, III.$\cC$.4.17 in \cite{BH:NPC:99}]\label{th:developability_of_nonpositively_curved_complexes}
	A non-positively curved complex of groups is developable.
\end{Theorem}

\noindent We will use the following criterion later on.

\begin{Definition}[I.7.38 in \cite{BH:NPC:99}]\label{def:geometric_link}
    For each point $x$ in a piecewise Euclidean polyhedral complex $X$, consider the set $\lk(x,X)$ of germs of $\sim$-equivalence classes of geodesic segments $[x,y]$ with $y\in \st(x)\setminus\{x\}$, $[x,y]\sim [x,y']$ if $[x,y]\subseteq [x,y']$ or vice versa.

    It can be endowed with the following metric: for two segments $[x,y]$, $[x,y']$ where $y$, $y'$ are contained in a common polyhedral cell, their distance is the comparison angle $\overline\angle_x(y,y')$ in the corresponding Euclidean comparison triangle. The metric on $\lk(x,X)$ is then the quotient metric associated to the projection $X\rightarrow K$.

    The set $\lk(x,X)$ has then the structure of a piecewise spherical polyhedral complex, called the \emph{geometric link at $X$}. If $\varepsilon>0$ is less than the distance from $x$ to the union of faces of $X$ that do not contain $x$, then the $\varepsilon$-neighbourhood of $x$ in $X$ is isometric to the $\varepsilon$-neighbourhood of the vertex $x$ in the Euclidean cone $C(\lk(x,X))$.
\end{Definition}

\begin{Proposition}[III.$\cC$.4.18 in \cite{BH:NPC:99}]\label{prop:geometric_link_nonpos_curved}
	If the geometric link $\lk(\tilde v,\st(\tilde v))$ is CAT(1) for every vertex $v\in V(\cY)$, then the complex of groups is non-positively curved and hence developable.
\end{Proposition}

\section{Two-dimensional affine buildings}\label{sec:buildings}

In this section, we will briefly introduce spherical and affine buildings. For general introductions to the topic, see \cite{Bro:Bdg:89} or \cite{Ron:LoB:89}. Furthermore, we will cite a criterion by Charney and Lytchak to recognise affine buildings by their metric properties.

\subsection{Definitions}

\begin{Definition}
	A \emph{Coxeter group $W$ of rank $r$} is a finitely presented group with the presentation
	\[
	W = \langle s_1,\ldots,s_r \,|\, (s_is_j)^{m(i,j)}=1 \rangle,
	\]
	where the \emph{Coxeter matrix $(m(i,j))_{i,j}$} is symmetric with value 1 on the diagonal and values in $\{2,3,\ldots,\infty\}$ else. We will denote the set of generators by $S=\{s_1,\ldots,s_r\}$.

	The Coxeter group $W$ is called \emph{irreducible} if there is no decomposition $S=S_1\sqcup S_2$ such that $W=\langle S_1\rangle\times \langle S_2\rangle$.
\end{Definition}

\begin{Definition}
    As usual, we associate to each Coxeter matrix a labelled graph with vertex set $S$, the \emph{Coxeter diagram}. Two vertices $s_i,s_j$ are connected by an edge if $m(i,j)>2$ and the edges are labelled with $m(i,j)$ if $m(i,j)>3$. In this paper, we will only consider Coxeter groups of types $\tilde A_2$, $\tilde C_2$ and $I_2(n)$ with the following Coxeter diagrams
	
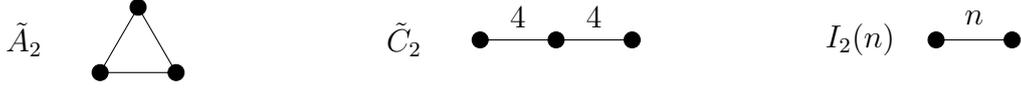
\begin{center}\begin{tikzpicture}
		\tikzstyle{every node}=[draw,circle,fill,inner sep=0pt,minimum width=6pt];
		\draw (0,-.433) node {} -- (1,-.433) node {} -- (.5,.433) node {}  -- cycle;
		\draw (5,0) node {} -- (6,0) node {} -- (7,0) node {} ;
		\draw (11,0) node {} -- (12,0) node {};
		\tikzstyle{every node}=[];
		\node (a2) at (-1,0) {$\tilde A_2$};
		\node (c2) at (4,0) {$\tilde C_2$};
        \node (i2) at (10,0) {$I_2(n)$};
		\node (1) at (5.5,.3) {$4$};
		\node (1) at (6.5,.3) {$4$};
		\node (1) at (11.5,.3) {$n$};
	\end{tikzpicture}\end{center}
\end{Definition}

\noindent Coxeter complexes are the basic ingredients we require for the definition of buildings.

\begin{Definition}
	The set of all cosets of special subgroups
	\[
		\Sigma(W,S)=\{w\langle T\rangle: w\in W, T\subseteq S\},
	\]
	partially ordered by reverse inclusion forms a simplicial complex, the \emph{Coxeter complex} associated to $(W,S)$. Top-dimensional simplices of $\Sigma(W,S)$ are called \emph{chambers}, codimension-1-simplices are called \emph{panels}.

	\begin{itemize}
		\item If $W$ is finite, in particular in the case $I_2(n)$, then there is a natural geometric realisation of $\Sigma(W,S)$ as a triangulated unit sphere $S^{r-1}$ and there is hence a piecewise spherical metric on $|\Sigma(W,S)|$. In this case $W$ is called \emph{spherical}.
        \item If the Coxeter diagram is \emph{affine}, in particular in the cases $\tilde A_2$ and $\tilde C_2$ (see \cite{Wei:SAB:09} for a complete list of affine diagrams), there is a natural geometric realisation of $\Sigma(W,S)$ as a triangulated Euclidean space $\R^{r-1}$. There is hence a piecewise Euclidean metric on $|\Sigma(W,S)|$, and $W$ is said to be \emph{affine} or irreducible Euclidean.
	\end{itemize}
\end{Definition}

\begin{Definition}
	A simplicial complex $X$ together with a collection of subcomplexes $\cA$ called \emph{apartments} is a \emph{building} if
	\begin{itemize}
		\item Every subcomplex $A\in\cA$ is a Coxeter complex.
		\item Every two simplices of $X$ are contained in a common apartment.
		\item For any two apartments $A_1,A_2\in \cA$ containing a common chamber, there is an isomorphism $A_1\rightarrow A_2$ fixing $A_1\cap A_2$ pointwise.
	\end{itemize}
	The building $X$ is called \emph{thick} if every panel is contained in at least three chambers. The axioms force all apartments to be isomorphic, there is, in particular a unique Coxeter group $W$ associated to $X$. Again, if $W$ is spherical or affine, we say that $X$ is spherical or affine.
\end{Definition}

\begin{Remarks}
    A building of type $I_2(n)$ is called a \emph{generalised $n$-gon}. These can also be described in a more explicit fashion, see Definition \ref{def:polygon}.

	There is a colouring of the vertex set of $X$ by the elements of $S$. The link of every vertex of $X$ is a building of the type obtained by removing the vertex colour from the generating set $S$. For the buildings of type $\tilde A_2$ and $\tilde C_2$, this implies that vertex links are generalised polygons.
\end{Remarks}

\subsection{Buildings as metric spaces}

\noindent The natural metric on $|\Sigma(W,S)|$ extends to a natural metric on $|X|$. If $W$ is spherical (affine), then $|X|$ is hence a piecewise spherical (Euclidean) simplicial complex.

\begin{Theorem}[Proposition 2.3 in \cite{CL:MC:01}]\label{th:buildings_are_cat}
	The canonical realisation of a spherical building is CAT(1) and has diameter $\pi$. The canonical realisation of an affine building is CAT(0). For every vertex $v$ in an affine building $X$, the geometric link $\lk(x,|X|)$ is a spherical building and the induced metric coincides with the natural building metric.
\end{Theorem}

It turns out that spherical and affine buildings can be recognised by their metric properties. This is the content of a paper by Charney and Lytchak, see \cite{CL:MC:01}. We will only require the following result, which is a special case of a result in this paper.

\begin{Theorem}[Theorem 7.3 in \cite{CL:MC:01}]\label{th:recognition}
	Let $X$ be a connected two-dimensional piecewise Euclidean polyhedral complex satisfying the following conditions:
	\begin{itemize}
		\item $X$ is CAT(0).
		\item Every 1-cell is contained in at least three 2-cells.
		\item The geometric link of every vertex is connected and has diameter $\pi$.
	\end{itemize}
	Then $X$ is an irreducible two-dimensional Euclidean building or a product of two trees.
\end{Theorem}

\subsection{Exotic buildings}\label{subsec:exotic}

The classification of affine buildings (due to Bruhat and Tits, see the book by Weiss \cite{Wei:SAB:09}) asserts that all affine buildings in dimension 3 and higher come from algebraic data such as an algebraic group over a local field.
In dimension 2, there is a free construction of other, exotic buildings by Ronan in \cite{Ron:CBR:86} which is very general but provides no control over the automorphism group. This construction in particular implies that there cannot be any hope to classify these buildings, because any classification would have to include the classification of all projective planes and generalised quadrangles.

In a locally finite affine building of type $\tilde A_2$, the number $q+1$ of chambers containing a common panel is constant throughout the building. For every prime power $q=p^e$, there are two classical buildings associated to the projective special linear groups over either a finite extension of the $p$-adic numbers $\Q_p$ or over the Laurent series field $\F_q(\!(t)\!)$. However, there are many more so-called $\emph{exotic buildings}$ of type $\tilde A_2$, some with very large automorphism groups. In \cite{HvM:NTB:90}, for example, van Maldeghem constructs exotic buildings of type $\tilde A_2$ with a vertex-transitive automorphism group.

The situation is similar in the case of buildings of type $\tilde C_2$. Here, there are more different types of classical buildings, but again there are uncountably many exotic buildings of type $\tilde C_2$. In \cite{Kan:GPS:86}, Kantor constructs exotic buildings of type $\tilde C_2$ with a large automorphism group.

\section{Singer polygons}\label{sec:singer_polygons}

We will construct complexes of groups whose universal covers are two-dimensional buildings. We will do this by prescribing the vertex links, which are generalised polygons. The main ingredients for our constructions are generalised polygons admitting \emph{Singer groups}.

\begin{Definition}\label{def:polygon}
    A \emph{generalised $n$-gon} is an incidence structure $\cI=(\cP,\caL,\cF)$ whose associated incidence graph has diameter $n$ and girth $2n$. The associated simplicial complex with vertex set $\cP\sqcup \caL$ and with 1-simplices $\cF$ is a spherical building of type $I_2(n)$ in the sense of section \ref{sec:buildings}.
\end{Definition}

\begin{Definition}
    A \emph{Singer polygon} is a generalised polygon which admits an automorphism group acting regularly, that is transitively and freely, on its points, a so-called \emph{Singer group}.
\end{Definition}

\noindent We will investigate examples for Singer polygons in the following sections.

\subsection{Bipartite graphs}

A generalised 2-gon is just a complete bipartite graph.

\begin{Lemma}
	Any group of order $k$ acts point- and line-regularly on a bipartite graph of order $(k,k$).
\end{Lemma}

\begin{Proof}
	Identify the sets of points and lines of the bipartite graph each with the elements of the group and consider the left regular representation on both sets.
\end{Proof}

\subsection{Projective planes}\label{subsec:proj_planes}

A generalised triangle is a projective plane. A projective plane is \emph{classical} if it is the usual projective plane over a division ring. The following classical result gives rise to the name of Singer groups.

\begin{Theorem}[Singer, \cite{Si:FPG:38}]
    Every finite classical projective plane admits a cyclic Singer group. A generator of such a group is called a \emph{Singer cycle}.
\end{Theorem}

\noindent It is conjectured that, conversely, every finite projective plane even admitting a group acting transitively on points is already classical. This seems to be difficult to prove, see \cite{Gil:TPP:07} for a recent contribution. It is not even known whether a projective plane with an abelian Singer group must be classical. On the other hand, all Singer groups on finite classical projective planes are classified by Ellers and Karzel:

\begin{Theorem}[Ellers-Karzel, 1.4.17 in \cite{De:FG:68}]
    Every finite classical projective plane of order $p^e$ admits a Singer group $S$ as follows:

    There is a divisor $n$ of $3e$ such that $n(p^e-1)$ divides $p^{3e}-1$ and $S$ can be presented in the following way:
    \[
    S = \langle \gamma, \delta \,|\, \gamma^s=1, \delta^n=\gamma^t, \delta\gamma\delta^{-1}= \gamma^{p^{3e/n}} \rangle,
    \]
    where $s = \frac{p^{3e}-1}{n(p^e-1)}$ and $t = \frac{p^{3e}-1}{n(p^{3e/n}-1)}$.

    Conversely, every group satisfying these relations operates point-regularly on a classical projective plane. In addition, the group $S$ is contained in $\PGl_3(\F_{p^e})$ if and only if $n\in\{1,3\}$.
\end{Theorem}

Finally, we require the following simple fact about arbitrary finite projective planes.

\begin{Proposition}[4.2.7 in \cite{De:FG:68}]\label{prop:regular_on_lines}
	Every Singer group on a finite projective plane also acts regularly on lines.
\end{Proposition}

\subsection{Generalised quadrangles}\label{subsec:GQs}

There is no classification of all Singer quadrangles known to the author. However, there is a list of all Singer quadrangles among all known finite generalised quadrangles, see the book \cite{STW:SQ:09}. All known Singer quadrangles except one arise by Payne derivation from translation generalised quadrangles and symplectic quadrangles.

Since we require rather explicit descriptions of the quadrangles, we consider only a large class of examples, the so-called \emph{slanted symplectic quadrangles}. For a detailed investigation of these quadrangles, see \cite{GJS:SSQ:94} and \cite{Str:PSQ:03}.

\begin{Definition}
     Fix a prime power $q>2$. Let $V=\langle e_{-2},e_{-1},e_1,e_2\rangle$ be a four-dimensional vector space over the finite field $\F_q$ endowed with the symplectic form $h$ satisfying
    \[
        h(e_i,e_j) =\begin{cases}
            -1 & i+j=0, i>j \\
            1 & i+j=0, i<j \\
            0 & \text{else.}
        \end{cases}
    \]
    The polar space consisting of all nontrivial totally isotropic subspaces of $V$, that is, subspaces $U$ satisfying $U\subseteq U^\perp$, forms a generalised quadrangle, called the \emph{symplectic quadrangle $W(q)$}. More specifically, set
    \begin{align*}
        \cP(W(q)) &= \{ U\subseteq V : \dim(U)=1, U\subseteq U^\perp \}, \\
        \caL(W(q)) &= \{ W\subseteq V : \dim(W)=2, W\subseteq W^\perp \}, \\
        \cF(W(q)) &= \{ (U,W)  \in \cP(W(q)) \times \caL(W(q)) : U\subseteq W \}.
    \end{align*}
\end{Definition}

We denote the full symplectic group by $\Sp_4(q)$, which is the group of all linear maps preserving $h$, and we fix the following notation
\[
    x(a) =\begin{pmatrix}
        1 & a & 0 & 0 \\
        0 & 1 & 0 & 0 \\
        0 & 0 & 1 & -a \\
        0 & 0 & 0 & 1
    \end{pmatrix},\qquad y(b) =\begin{pmatrix}
        1 & 0 & b & 0 \\
        0 & 1 & 0 & b \\
        0 & 0 & 1 & 0 \\
        0 & 0 & 0 & 1
    \end{pmatrix},\qquad z(c) =\begin{pmatrix}
        1 & 0 & 0 & c \\
        0 & 1 & 0 & 0 \\
        0 & 0 & 1 & 0 \\
        0 & 0 & 0 & 1
    \end{pmatrix}
\]
for specific elements in $\Sp_4(q)$, where $a,b,c\in\F_q$. We abbreviate $x=x(1)$, $y=y(1)$ and $z=z(2)$. It is not hard to see that $[x,y]=z$.

\begin{Lemma}\label{lemma:heisenberg_group}
    The subgroup $E=\langle x(a),y(b),z(c) : a,b,c\in\F_q \rangle\leq \Sp_4(q)$ acts regularly on all points not collinear with the point $p_0=\F_q(1,0,0,0)^T$.
    \begin{itemize}
	    \item If $q$ is odd, then $E$ is a three-dimensional Heisenberg group over $\F_q$. If $q$ is prime, we have the following simple presentation:
    \[
        E = \langle x,y \,|\, z=xyx^{-1}y^{-1}, x^q=y^q=z^q=1, xz=zx, yz=zy \rangle
    \]
        \item If $q$ is even, then $E\cong \F_q^3$, which is hence an elementary abelian 2-group.
    \end{itemize}
\end{Lemma}

\begin{Proof}
    It is an easy calculation to show that all points not collinear with $p_0$ have the form $\F_q(x,y,z,1)^T$. Simple matrix calculations then imply that $E$ acts regularly on these points.

    If $q$ is odd, then $E$ is a three-dimensional Heisenberg group by \cite[4.3]{STW:SQ:09}. It is a simple calculation to verify the presentation in case that $q$ is prime. By calculating three commutators, one can see that $E$ is isomorphic to $\F_q^3$ if $q$ is even.
\end{Proof}

\paragraph{Payne derivation} We will give a definition of \emph{Payne derivation} as described in \cite{STW:SQ:09}. The same procedure is called \emph{slanting} in \cite{GJS:SSQ:94}. We write $p\sim r$ for collinear points $p$ and $r$. Then we define
\begin{align*}
    p^\perp & = \{ r\in \cP(W(q)) : p \sim r \} \\
    \{p,p'\}^{\perp\perp} &= \{ r\in\cP(W(q)) : r \in s^\perp \text{ for all } s\in p^\perp \cap p'^\perp \}.
\end{align*}
\begin{Definition}
    The \emph{slanted symplectic quadrangle $W(q)^\Diamond$} is given as follows:
\begin{itemize}
    \item The points of $W(q)^\Diamond$ are all points of $W(q)$ not collinear with $p_0$.
    \item The lines of $W(q)^\Diamond$ are all lines of $W(q)$ not meeting $p_0$ as well as the sets
	    \[\{p_0,r\}^{\perp\perp} \setminus \{p_0\}\] for all points $r\in W(q)^\Diamond$.
\end{itemize}
\end{Definition}

\noindent The following characterisation will be used later on to construct lattices in buildings of type $\tilde C_2$.

\begin{Theorem}\label{th:slanted_sympl_quadrangle}
	The slanted symplectic quadrangle $W(q)^\Diamond$ is a Singer quadrangle of order $(q-1,q+1)$ with Singer group $E$. A set of representatives $L$ for the $E$-action on lines is given by all lines through the point $p_1=(0,0,0,1)^T$:
	\[
    L = \bigl\{ l_{[a:b]} = \bigl(\F_q(0,0,0,1)^T + \F_q (0,b,a,0)^T\bigr) \,\big|\, [a:b]\in \bP\F_q^2 \bigr\} \sqcup \bigl\{l_0 = \{p_0,p_1\}^{\perp\perp}\setminus\{p_0\}\bigr\}
	\]
	The line stabilisers have the form
	\[
	E_{l_{[a:b]}}=\biggl\{\begin{pmatrix}
		1 & fa & fb & 0 \\
		0 & 1 & 0 & fb \\
		0 & 0 & 1 & -fa \\
		0 & 0 & 0 & 1
    \end{pmatrix}: f\in \F_q \biggr\}\quad\text{and}\quad E_{l_0} = \{ z(f) : f\in\F_q \}.
	\]
	In particular, all line stabilisers are isomorphic to the additive group of $\F_q$. If $q$ is an odd prime, the line stabilisers have the form
	\[
    E_{l_{[a:b]}} = \langle x^ay^bz^{-\tfrac{1}{2}ab} \rangle\quad\text{and}\quad E_{l_0} = \langle z \rangle,
	\]
    where of course $2$ is invertible in $\F_q^\times$. If $q$ is even, we have
    \[
    E_{l_{[a:b]}} = \langle x(a)y(b) \rangle\quad\text{and}\quad E_{l_0} = \langle z(1) \rangle
    \]
    as $\F_q$-subspaces in $\F_q^3$.
\end{Theorem}

\begin{Proof}
	The first claim follows from Lemma \ref{lemma:heisenberg_group}. For the set of line representatives and the structure of the stabilisers, see \cite[Lemma 4.1]{Str:PSQ:03}. The last two statements are simple calculations.
\end{Proof}

\section{A local complex of groups construction}\label{sec:local_cplx_of_groups_construction}

Using Singer polygons, we will now construct a complex of groups whose local development at one vertex is a generalised polygon. With this local piece of data, we will later construct complexes of groups whose local developments are buildings.

\begin{Definition}
	Associated to any generalised $n$-gon $\cI=(\cP,\caL,\cF)$, there is a scwol $\cZ(\cI)$ with vertex set
	\begin{align*}
		V(\cZ(\cI)) &= \cP\sqcup\caL\sqcup\cF \\
		\intertext{and with edges, suggestively written as follows,}
		E(\cZ(\cI)) &= \{ p \leftarrow f : p\in\cP, f\in \cF, p\in f \} \sqcup \{ l\leftarrow f: l\in\caL, f\in \cF,l\in f\},
\end{align*}
where we mean $i(x\leftarrow y)=y$ and $t(x\leftarrow y)=x$, of course.
\end{Definition}

\noindent The geometric realisation of this scwol is isomorphic as a simplicial complex to the barycentric subdivision of the geometric realisation $|\cI|$ of the generalised polygon.

\paragraph{Notation} Let $\cI=(\cP,\caL,\cF)$ be a Singer polygon with Singer group $S$. Fix a point $p\in \cP$ and a set $L$ of line representatives for the $S$-action on $\caL$. For all lines $l\in L$, we consider the sets
\[
	D_l = \{ d\in S: (p,dl)\in\cF \}.
\]
The set $D = \bigcup_{l\in L} D_l$ is called a \emph{difference set}. Furthermore, let $F_l\coloneq \{ (p,dl) : d\in D_l\}$ for all $l\in L$. This is a partition of the set of all flags containing $p$, which we denote by
\[
F \coloneq \{f\in\cF : p\in f\} = \coprod_{l\in L}F_l = \coprod_{l\in L} \{ (p,dl):d\in D_l\}.
\]
In particular, the set $F$ is a set of representatives of flags for the $S$-action on $\cF$.

\paragraph{Construction} Let a complex of groups $G(\cY)$ be given. Let $v\in V(\cY)$ be a minimal element (that is: no edges start in $v$), and assume that the closed star of $v$ in $\cY$, denoted by $\cY(v)$, has the following structure:
\begin{align*}
	V(\cY(v)) \cong&\,\{v \} \sqcup \{p\} \sqcup L \sqcup F. \\
	\intertext{Again, we will write the edges in the following suggestive fashion:}
	E(\cY(v))\cong&\,\{ v\leftarrow p \} \sqcup \{ v\leftarrow l : l\in L \} \sqcup \{ v\leftarrow f : f  \in F \} \\
	& \quad\sqcup \{ p\leftarrow f : f\in F \} \sqcup \coprod_{l\in L}\Bigl\{ l\leftarrow f: f\in F_l\Bigr\}.
\end{align*}

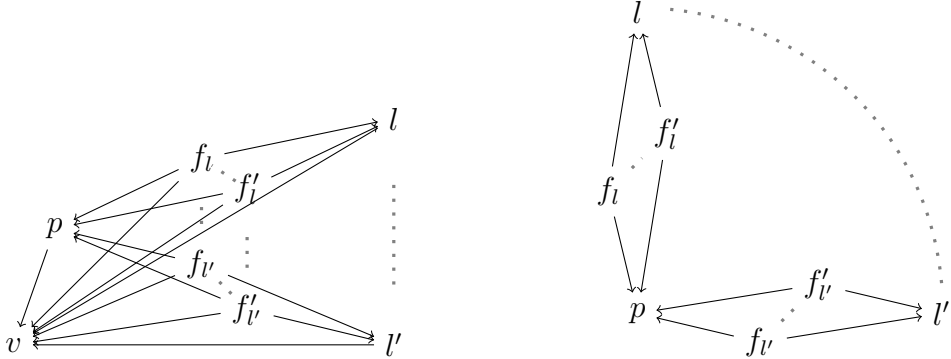
\begin{figure}[hbt]
\centering
    \begin{tikzpicture}
	\node (v) at (0,0,0) {$v$};
	\node (p) at (-1,0,-4) {$p$} edge[->] (v);
	\node (l1) at (5,3,0) {$l$} edge[->] (v);
	\node (l2) at (5,0,0) {$l'$} edge[->] (v) edge[color=gray, loosely dotted,very thick, shorten <=.5cm, shorten >=.5cm] (l1);
	\node (f1a) at (1.7,1.7,-2) {$f_{l}$} edge[->] (v) edge[->] (p) edge[->] (l1);
	\node (f1b) at (2.3,1.3,-2) {$f'_{l}$} edge[->] (v) edge[->] (p) edge[->] (l1) edge[color=gray, loosely dotted,very thick, shorten <=.2cm, shorten >=.2cm] (f1a);
	\node (f2a) at (1.7,0.3,-2) {$f_{l'}$} edge[->] (v) edge[->] (p) edge[->] (l2) edge[color=gray, loosely dotted, very thick, shorten <=.2cm, shorten >=.2cm] (f1a);
	\node (f2b) at (2.3,-0.3,-2) {$f'_{l'}$} edge[->] (v) edge[->] (p) edge[->] (l2) edge[color=gray, loosely dotted,very thick, shorten <=.2cm, shorten >=.2cm] (f2a) edge[color=gray, loosely dotted, very thick, shorten <=.2cm, shorten >=.2cm] (f1b);
\end{tikzpicture} \hspace{2cm}
\begin{tikzpicture}[scale=2]
	\node (p) at (0,0,0) {$p$};
	\node (l1) at (2,0,0) {$l'$};
	\node (l2) at (0,2,0) {$l$};
	\node (f1a) at (1,0,.5) {$f_{l'}$} edge[->] (p) edge[->] (l1);
    \node (f1b) at (1,0,-.5) {$f'_{l'}$} edge[->] (p) edge[->] (l1) edge[color=gray,loosely dotted,very thick, shorten <=.2cm, shorten >=.2cm] (f1a);
	\node (f2a) at (0,1,.5) {$f_l$} edge[->] (p) edge[->] (l2);
    \node (f2b) at (0,1,-.5) {$f'_{l}$} edge[->] (p) edge[->] (l2) edge[color=gray,loosely dotted,very thick, shorten <=.2cm, shorten >=.2cm] (f2a);
	\draw[color=gray,loosely dotted,very thick] (2,.2,0) arc (5:85:2);
\end{tikzpicture} \caption{The scwol $\cY(v)$ on the left, the link $\Lk_v(\cY)$ on the right.}\label{fig:local_development}
\end{figure}

The scwol $\cY(v)$ can be visualised as in Figure \ref{fig:local_development}. The left image shows the full complex $\cY(v)$, while the right picture illustrates the upper link of the vertex $v$.

In addition, assume now that the vertex groups in $G(\cY)$ are as follows
\[
G_v \cong S,\qquad G_p=G_{f}=\{1\},\qquad G_l \cong S_l \qquad \forall l\in L, \forall f\in F,
\]
that the monomorphisms are the obvious inclusions and that the twist elements are
\[
g_{v\leftarrow p\leftarrow (p,dl)} = 1, \qquad g_{v\leftarrow l \leftarrow (p,dl) } = d^{-1}\qquad\forall l\in L, d\in D_l.
\]
The significance of this construction lies in the following

\begin{Proposition}\label{prop:local_development_v}
	The local development $\cY(\tilde v)$ isomorphic to the cone over the scwol $\cZ(\cI)$.
\end{Proposition}

\begin{Proof}
	Since $v$ is a minimal element in $\cY$, by Definition \ref{def:local_development}, we have $\cY(\tilde v) = \{\tilde v\} * \Lk_{\tilde v}(\cY)$. By the construction of the scwol $\Lk_{\tilde v}(\cY)$ in Definition \ref{def:upper_link}, we have
	\begin{align*}
		V(\Lk_{\tilde v}(\cY)) =& \,\bigl\{ (g\psi_a(G_{i(a)}),a) : a\in E(\cY), t(a)=v, g\psi_a(G_{i(a)})\in G_v/\psi_a(G_{i(a)}) \bigr\} \\
		=&\, \{ ( \{g\}, v \leftarrow p ) : g\in S\} \sqcup \{ (g S_l, v\leftarrow l) : g\in S, l\in L \}\\
		&\,\quad \sqcup \{ (\{g\}, v\leftarrow f) : g\in S, f\in F \},
	\end{align*}
	which is in bijection to $V(\cZ(\cI))=\cP \sqcup \caL \sqcup \cF$ via $(\{g\},v\leftarrow p)\mapsto gp$, $(gS_l,v\leftarrow l)\mapsto gl$ and $(\{g\}, v\leftarrow f)\mapsto gf$. For the edge set, we obtain
	\begin{align*}
		E(\Lk_{\tilde v}(\cY)) =& \,\bigl\{ (g\psi_{ab}(G_{i(b)}),a,b) : (a,b) \in E^{(2)}(\cY), t(a)=v, g\psi_{ab}(G_{i(b)}) \in G_v/\psi_{ab}(G_{i(b)}) \bigr\} \\
		=&\, \{ ( \{g\}, v\leftarrow p \leftarrow f) : g\in G, f\in F \} \\
		&\,\quad \sqcup \{ ( \{g\}, v\leftarrow l \leftarrow f) : g\in G, l\in L, f\in F_l \}
	\end{align*}
	which is in bijection to the set of edges $E(\cZ(\cI))$ via
	\begin{align*}
        ( \{g\}, v\leftarrow p \leftarrow f ) &\mapsto (gp \leftarrow gf ) \\
        ( \{g\}, v\leftarrow l \leftarrow (p,dl) ) &\mapsto (gdl \leftarrow (gp,gdl) ).
	\end{align*}
	The twist elements precisely guarantee that these bijections commute with the maps $i$ and $t$, respectively. We have
	\[\begin{array}{rclcrcl}
		i( (\{g\},v\leftarrow p \leftarrow f)) &=& (\{g\}, v\leftarrow f) &\,\mapsto\,& gf &=& i(gp\leftarrow gf) \\
		t( (\{g\},v\leftarrow p \leftarrow f)) &=& (\{g\}, v\leftarrow p) &\,\mapsto\,& gp &=& t(gp\leftarrow gf) \\
		i( (\{g\},v\leftarrow l \leftarrow (p,dl))) &=& (\{g\}, v\leftarrow (p,dl)) &\,\mapsto\,& (gp,gdl) &=& i(gdl\leftarrow (gp,gdl)) \\
		t( (\{g\},v\leftarrow l \leftarrow (p,dl))) &=& (gdS_l, v\leftarrow l) &\,\mapsto\,& gdl &=& t(gdl\leftarrow (gp,gdl)),
	\end{array}\]
	where we use Definition \ref{def:upper_link} for the maps $i$ and $t$. So $\Lk_{\tilde v}(\cY)\cong\cZ(\cI)$.
\end{Proof}

\noindent Note that the special situation is always simpler.
    \begin{description}
        \item[Bipartite graphs] If a group $S$ acts regularly on points and lines of a complete bipartite graph of order $(k,k)$, there is only one line orbit, hence $L=\{l\}$ and $S_l=\{1\}$. Since the bipartite graph is complete, we have $|F|=|F_l|=k$ and $D_l=S$.
        \item[Projective planes] If a group $S$ is a Singer group on a projective plane, there is only one line orbit by Proposition \ref{prop:regular_on_lines}, hence $L=\{l\}$ and $S_l=\{1\}$. In particular, there is only one difference set $D\coloneq D_l$.
        \item[Generalised quadrangles] For slanted symplectic quadrangles of order $(q-1,q+1)$, there are $q$ line orbits, hence $|L|=|F|=q$. For simplicity, we choose the set $L$ to be the set of all lines incident to $p$. Then we have $D_l=\{1\}$ for all $l\in L$.
    \end{description}

\noindent We endow the geometric realisation $|\cY|$ with a piecewise Euclidean metric such that the angles between $v\leftarrow p$ and $v\leftarrow f$ as well as the angles between $v\leftarrow f$ and $v\leftarrow l$ are $\pi/2n$ for all lines $l\in L$ and all flags $f\in F$.

\begin{Proposition}\label{prop:key}
    Then the geometric link $\lk(\tilde v,\st(\tilde v))$ is isometric to the barycentric subdivision of the generalised $n$-gon $|\cI|$ with its standard metric. In particular, it is a connected CAT(1) space with diameter $\pi$.
\end{Proposition}

\begin{Proof}
    By Proposition \ref{prop:local_development_v}, the local development $\cY(\tilde v)$ is isomorphic to the cone over the scwol $\cZ(\cI)$. By Definition \ref{def:geometric_link}, the polyhedral complex structure on $\lk(\tilde v,\cY(\tilde v))$ is then the barycentric subdivision of $\cI$.
    By construction, the angles at the vertex $\tilde v$ are $\pi/2n$. Since the edge length in the geometric link $\lk(\tilde v,\st(\tilde v))$ is given by the angles, the geometric link is hence isometric to the barycentric subdivision of $|\cI|$ with the natural metric of a generalised polygon. The CAT(1) condition follows from Theorem \ref{th:buildings_are_cat}.

\end{Proof}

\section{Lattices in buildings of type \texorpdfstring{$\tilde A_2$}{\textasciitilde A2}}\label{sec:a2}

In this section, we will give a construction of small complexes of groups whose universal covers are affine buildings of type $\tilde A_2$. We will first state a general construction. Later, we specialise to cyclic Singer groups to obtain lattices with very simple presentations. For these special lattices, we investigate the structure of spheres of radius 2 in the corresponding buildings.

\subsection{A general construction of panel-regular lattices in buildings of type \texorpdfstring{$\tilde A_2$}{\textasciitilde A2}}\label{subsec:general_a2_construction}

For the whole construction, we fix three finite Singer projective planes $\cI_1$, $\cI_2$ and $\cI_3$ of order $q$ with three (possibly isomorphic) Singer groups $S_1$, $S_2$ and $S_3$. Note that the projective planes need not be classical, even though this is likely in view of what we have said in section \ref{subsec:proj_planes}. By Proposition \ref{prop:regular_on_lines}, these Singer groups act regularly on lines as well. As in section \ref{sec:local_cplx_of_groups_construction}, we obtain a difference set for each of these groups by choosing a point and a line in the corresponding projective plane. We denote these three difference sets by $D_1$, $D_2$ and $D_3$.

Write $J=\{0,1,\dots,q\}$ and choose bijections $d_\alpha : J\rightarrow D_\alpha$ for $\alpha\in \{1,2,3\}$, which we call \emph{ordered difference sets}.

\paragraph{Construction} Associated to this piece of data, we consider the complex of groups $G(\cY)$ over the scwol $\cY$ with vertices
\begin{align*}
       V(\cY) &\coloneq \{ v_1, v_2, v_3\} \sqcup \{e_1,e_2,e_3\} \sqcup \{ f_j : j \in J \} \\
       \intertext{and edges}
       E(\cY) &\coloneq \{ v_\alpha \leftarrow e_\beta: \alpha\neq\beta\} \sqcup \{ v_\alpha \leftarrow f_j : j\in J\} \sqcup \{e_\beta \leftarrow f_j : j\in J\}.
\end{align*}

\noindent Figure \ref{fig:cYfor2} illustrates this scwol for $q=2$.

\begin{figure}[hbt]
\centering
\begin{tikzpicture}[scale=4,font=\small]
       \node (p1) at (0,0,0) {$v_1$};
       \node (p2) at (2,0,0) {$v_2$};
       \node (p3) at (1,0,.866) {$v_3$};
       \node (l3) at (1,0,0) {$e_3$} edge[->] (p1) edge[->] (p2);
       \node (l1) at (1.5,0,.433) {$e_1$} edge[->] (p2) edge[->] (p3);
       \node (l2) at (.5,0,.433) {$e_2$} edge[->] (p1) edge[->] (p3);
       \node (f1) at (1,.5,.433) {$f_1$} edge[->] (p1) edge[->] (p2) edge[->] (p3) edge[->] (l1) edge[->] (l2) edge[->] (l3);
       \node (f2) at (1,0,.433) {$f_2$} edge[->] (p1) edge[->] (p2) edge[->] (p3) edge[->] (l1) edge[->] (l2) edge[->] (l3);
       \node (f3) at (1,-.5,.433) {$f_3$} edge[->] (p1) edge[->] (p2) edge[->] (p3) edge[->] (l1) edge[->] (l2) edge[->] (l3);
\end{tikzpicture} \caption{The scwol $\cY$ for $q=2$}\label{fig:cYfor2}
\end{figure}

Choose the three vertex groups $G_{v_\alpha}$ to be isomorphic to $S_\alpha$ for $\alpha=1,2,3$. All other vertex groups are trivial. The twist elements are defined as follows:
\[
g_{v_\alpha\leftarrow e_\beta \leftarrow f_j} \coloneq\begin{cases}
       1 & \beta - \alpha \equiv_3 1 \ \\
       d_\alpha(j)^{-1}& \beta - \alpha \equiv_3 2.
\end{cases}
\]

\noindent We endow $|\cY|$ with a locally Euclidean metric as follows:

Let $\Delta$ be the geometric realisation of one triangle in the affine Coxeter complex of type $\tilde A_2$. For each $j\in J$, we map each subcomplex spanned by $\{v_1,v_2,v_3,e_1,e_2,e_3,f_j\}$ onto the barycentric subdivision of $\Delta$ in the obvious way and pull back the metric. We obtain a locally Euclidean metric on $|\cY|$. In particular, the angles at every vertex $v_\alpha$ are $\pi/3$.

\begin{Proposition}\label{prop:a2_is_developable}
       The complex of groups $G(\cY)$ is developable.
\end{Proposition}

\begin{Proof} We want to apply Proposition \ref{prop:geometric_link_nonpos_curved}, so we have to check that the geometric link of every vertex in its local development is CAT(1).
    
    For any vertex $v_\alpha$, the subcomplex $\cY(v_\alpha)$ obviously has the appropriate structure for the construction in Section \ref{sec:local_cplx_of_groups_construction}. By construction, the induced complex of groups has the right form to apply Proposition \ref{prop:local_development_v}, and we obtain that $\cY(\tilde v_\alpha)$ is isomorphic to the cone over the scwol $\cZ(\cI_\alpha)$. Then, by Proposition \ref{prop:key}, the geometric link $\Lk(\tilde v_\alpha, \st(\tilde v_\alpha))$ is CAT(1).

       For any vertex $e_\beta$, by Definition \ref{def:local_development} the local development has the form
       \[
       \St(\tilde e_\beta) = \lvert\Lk^{e_\beta}(\cY) \ast\, \{\tilde e_\beta\} \ast\, \Lk_{\tilde e_\beta}(\cY)\rvert.
       \]
       The lower link $\Lk^{e_\beta}(\cY)$ is isomorphic to the scwol consisting of just two vertices $\{v_\alpha: \alpha\neq \beta\}$, the upper link of the local development $\Lk_{\tilde e_\beta}(\cY)$ is isomorphic to a scwol consisting of the vertices $\{f_j: j\in J\}$. The geometric realisation of the local development is hence isometric to $q+1$ triangles joined along one edge. In particular, the geometric link is CAT(1).

       The upper link of the local development $\Lk_{\tilde f_j}(\cY)$ is trivial. The lower link $\Lk^{f_j}(\cY)$ consists of the sub-scwol of $\cY$ spanned by the vertices $\{v_1,v_2,v_3,e_1,e_2,e_3\}$. The local development $\cY(\tilde f_j)$ is hence just one triangle, so the geometric link is also CAT(1). By Proposition \ref{prop:geometric_link_nonpos_curved}, the complex $G(\cY)$ is then non-positively curved and hence developable by Theorem \ref{th:developability_of_nonpositively_curved_complexes}.
\end{Proof}

\begin{Proposition}\label{prop:fundamental_group_a2}
       The fundamental group $\Gamma$ of the complex $G(\cY)$ admits the following presentation:
    \[
    \Gamma = \Bigl\langle S_1,S_2,S_3 \,\Big| \begin{array}{c}\text{ all relations in the groups }S_1,S_2,S_3,\\ d_1(j)d_2(j)d_3(j) = d_1(j')d_2(j')d_3(j') \quad\forall j,j'\in J\end{array}\Bigr\rangle
    \]
\end{Proposition}

\begin{Remark}
    We show in Theorem \ref{th:a2_classification} that $D_\alpha$ and $d_\alpha$ can always be chosen such that $d_\alpha(0)=1$, which simplifies this presentation even further.
\end{Remark}

\begin{Proof}
       Consider the following maximal spanning subtree $T$:
       \begin{align*}
               E(T)= \{ v_1\leftarrow e_2, v_1\leftarrow e_3, v_2\leftarrow e_3, v_2\leftarrow e_1, v_3\leftarrow e_1 \} \sqcup \{e_3\leftarrow f_j : j \in J \}.
       \end{align*}
    By taking a look at Definition \ref{def:fundamental_group}, we obtain a presentation of $\Gamma$ with generating sets $S_1$, $S_2$ and $S_3$ and generators for all edges of the scwol $\cY$ not contained in $T$. The relations imposed are first of all the relations from the groups $S_\alpha$. For the relations of the form $k_a k_b = g_{a,b}k_{ab}$, consider Figure \ref{fig:triangle_for_presentation} showing the $j$-th triangle in the complex of groups $G(\cY)$. There, black arrows indicate edges contained in $T$ and all other edges are drawn dotted. Group elements $g$ written on dotted edges $a$ indicate that $k_a=g$.

    \begin{figure}[hbt]
       \centering
       \begin{tikzpicture}[scale=6,font=\footnotesize]
               \node (v1) at (0,0) {$v_1$};
               \node (v3) at (1,0) {$v_3$};
               \node (v2) at (.5,.866) {$v_2$};
               \node (e2) at (.5,0) {$e_2$} edge[->] (v1) edge[->,dotted] node[below,color=gray]{$k_{v_3\leftarrow e_2}$} (v3);
               \node (e3) at (.25,.466) {$e_3$} edge[->] (v1) edge[->] (v2);
               \node (e1) at (.75,.466) {$e_1$} edge[->] (v2) edge[->] (v3);
               \node (fj) at (.5,.289) {$f_j$} edge[->] (e3) edge[->,dotted] node[color=gray,pos=0.7]{$d_1(j)$} (v1) edge[->,dotted] node[color=gray,above]{$1$} (v2) edge[->,dotted] node[pos=0.2,color=gray]{$d_2(j)^{-1}$} (v3) edge[->,dotted] node[color=gray]{$d_2(j)^{-1}$} (e1) edge[->,dotted] node[color=gray]{$d_1(j)$} (e2);
               \node (m12) at (.593,.540) {$d_2(j)^{-1}$};
               \node (m32) at (.417,.540) {$1$};
               \node (m31) at (.25,.252) {$d_1(j)^{-1}$};
               \node (m23) at (.75,.252) {$1$};
               \node (m21) at (.333,.096) {$1$};
               \node (m21) at (.666,.096) {$d_3(j)^{-1}$};
       \end{tikzpicture} \caption{One triangle in $G(\cY)$}\label{fig:triangle_for_presentation}
\end{figure}
All edge elements but $k_{v_3\leftarrow e_2}$ have already been replaced by elements in the groups $S_\alpha$. For $k_{v_3\leftarrow e_2}$, we obtain the relation $k_{v_3\leftarrow e_2}=d_1(j)d_2(j)d_3(j)$ for all $j\in J$. Since this edge is contained in all triangles for all $j\in J$, we obtain the additional relations for all pairs of triangles. All other relations from Definition \ref{def:fundamental_group} do not show up here, since almost all vertex groups are trivial.
\end{Proof}

\begin{Theorem}\label{th:a2_is_building}
The universal cover $\cX$ is a building of type $\tilde A_2$, where the vertex links are isomorphic to $\cI_1$, $\cI_2$ and $\cI_3$. The fundamental group $\Gamma=\pi_1(G(\cY))$ with presentation
    \[
    \Gamma = \Bigl\langle S_1,S_2,S_3 \,\Big| \begin{array}{c}\text{ all relations in the groups }S_1,S_2,S_3,\\ d_1(j)d_2(j)d_3(j) = d_1(j')d_2(j')d_3(j') \quad\forall j,j'\in J\end{array}\Bigr\rangle
    \]
is hence a uniform lattice in the full automorphism group of the building. The set of all chambers containing a fixed panel is a fundamental domain for the action.
\end{Theorem}

\begin{Proof}
        The universal cover $|\cX|$ is a simply connected space which is locally CAT(0). So by \cite[II.4.1]{BH:NPC:99}, the space $\lvert\cX\rvert$ is CAT(0). Consider the space $|\cX|$ as a polyhedral complex consisting only of simplices where each cell is the preimage of subscwols spanned by $\{v_1,v_2,v_3,e_1,e_2,e_3,f_j\}_{j\in J}$. The vertices of this polyhedral complex are the preimages of the vertices $v_{\alpha}$ of $\cX$. This is, a priori, not a simplicial complex since the intersection of two simplices might not be a simplex.

However, in view of Theorem \ref{th:recognition} it remains to see that $|\cX|$ is thick and that the geometric links of all vertices are connected and have diameter $\pi$. Thickness of $|\cX|$ is clear by construction. Again by Proposition \ref{prop:key}, the geometric link of every vertex is connected and of diameter $\pi$.

By Theorem \ref{th:recognition}, the geometric realisation $|\cX|$ is then either a two-dimensional affine building or a product of two trees. Of course, the universal cover cannot be a product of trees and it has to be of type $\tilde A_2$, since all vertex links are projective planes.
\end{Proof}

\begin{Remark}
    Unfortunately, there are many examples of buildings of type $\tilde A_2$. It is clear that $\cX$ has an automorphism group which is transitive on vertices of the same type, but there are known examples of non-classical buildings with vertex-transitive automorphism groups, see \cite{HvM:NTB:90}. However, some information can be obtained by considering spheres of radius 2 in $\cX$ and comparing these to corresponding spheres in the classical buildings of type $\tilde A_2$ associated to $\Q_p$ and $\F_p(\!(t)\!)$, respectively. We will explicitly calculate these spheres of radius 2 in section \ref{subsec:spheres}.
\end{Remark}

\paragraph{Construction} Since we know that $|\cX|$ is a building, it is a flag complex. By using the construction of the universal cover in section \ref{sec:complexes_of_groups}, we obtain a very explicit description of the building as the flag complex over the following graph $X$:
\begin{align*}
        V(X) &\coloneq \Gamma/S_1 \sqcup \Gamma/S_2 \sqcup \Gamma/S_3 \\
        E(X) &\coloneq \bigl\{ (gS_1,gS_2) : g\in\Gamma\bigr\} \sqcup
        \bigl\{ (gS_2,gS_3) : g\in\Gamma\bigr\}
		\\&\qquad\quad\sqcup \bigl\{ (gS_1,gk^{-1}S_3) : g\in\Gamma\bigr\},
\end{align*}
where $k^{-1}=d_1(0)d_2(0)d_3(0)$. Note again that we can choose $D_\alpha$ and $d_\alpha$ such that $k=1$ by Theorem \ref{th:a2_classification} to obtain a symmetric description.

\begin{Remark}
    Finally, it is not clear in which way the building $\cX$ and the lattice $\Gamma$ depend on the given data, namely on the Singer groups $S_\alpha$, on the difference sets $D_\alpha$ and on the ordered difference sets $d_\alpha$. In section \ref{subsec:cyclic_lattices}, we give examples where different orderings $d_\alpha$ lead to non-isomorphic lattices. We do not know whether the lattices arising from different orderings are all commensurable. In addition, it is not even clear whether this construction can lead to different buildings for fixed projective planes.
\end{Remark}

\noindent Several strong properties of these lattices can be obtained without knowing the exact building they act on.

\begin{Remark}
    By unpublished work of Shalom and Steger, lattices in arbitrary buildings of type $\tilde A_2$ follow a version of Margulis' normal subgroup theorem. In particular, the lattices $\Gamma$ are almost perfect, their first homology group is finite. For the special case of lattices generated by cyclic Singer groups, one can calculate group homology quite explicitly, as we show in section \ref{sec:group_homology}.
\end{Remark}

\noindent In addition, by using the so-called spectral criterion, one can see that these lattices have property (T).

\begin{Proposition}[Cartwright-M{\l}otkowski-Steger, \.Zuk]
	All cocompact lattices in buildings of type $\tilde A_2$ have property (T).
\end{Proposition}

\begin{Proof}
	See \cite{CMS:TA2:94} or \cite{Zuk:TGP:96}. An detailed description of the latter result can be found in \cite{BHV:KPT:08}.
\end{Proof}

\subsection{The Classification of panel-regular lattices and simple properties}\label{subsec:a2_classification}

In this section, we will classify all panel-regular lattices on locally finite buildings of type $\tilde A_2$. In this process, we will show that, actually, the examples from the last section cover all possible lattices, and that the presentations can even be simplified to a more symmetric form. We start with a simple observation.

\begin{Lemma}
    Let $X$ be a locally finite building of type $\tilde A_2$. Let $\Gamma\leq \Aut(X)$ be a lattice acting regularly on one type of panel. Then $\Gamma$ already acts regularly on all types of panels.
\end{Lemma}

\begin{Proof}
    Let $v,w,u$ be the vertices of a chamber in the building. Assume that $\Gamma$ acts regularly on the panels of the same type as the panel $(v,w)$. In particular, the vertex stabiliser $\Gamma_v$ is a group acting regularly on points (or lines) of the finite projective plane $\lk_X(v)$. By Proposition \ref{prop:regular_on_lines}, $\Gamma_v$ then also acts regularly on lines (or points) of $\lk_X(v)$, hence $\Gamma$ also acts regularly on the panels of the same type as $(v,u)$. By repeating the argument, we see that $\Gamma$ acts regularly on all types of panels.
\end{Proof}

In particular, the lattice $\Gamma$ acts transitively on vertices of the same type, which implies that vertex links of vertices of the same type are isomorphic. There are hence three finite projective planes of the same order, denoted by $\cI_1$, $\cI_2$ and $\cI_3$ associated to $X$. Pick a chamber in $X$ with vertices $\bar v_1, \bar v_2$ and $\bar v_3$. By the same argumentation as before, the three vertex stabilisers $\Gamma_{\bar v_1}$, $\Gamma_{\bar v_2}$ and $\Gamma_{\bar v_3}$ are then Singer groups on the vertex links and we denote these by $S_1$, $S_2$ and $S_3$. In the following, we write $J=\{0,1,\ldots,q\}$.

\begin{Theorem}\label{th:a2_classification}
    Let $X$ be a locally finite building of type $\tilde A_2$. We assume that there is a lattice $\Gamma\leq \Aut(X)$ acting regularly on one and hence on all types of panels. Then there are three Singer groups $S_1$, $S_2$ and $S_3$ associated to the three types of possible vertex links. Furthermore, there are three ordered difference sets $d_1$, $d_2$ and $d_3$ corresponding to these Singer groups satisfying $d_\alpha(0)=1\in S_\alpha$. The lattice $\Gamma$ admits the following presentation:
    \[
        \Gamma \cong \langle S_1,S_2,S_3 \,|\, \text{ all relations in the groups }S_1,S_2,S_3,\quad d_1(j)d_2(j)d_3(j) = 1 \quad\forall j \in J\rangle
    \]
\end{Theorem}

\begin{Proof}
    Consider the canonical scwol $\cX$ associated to the building. Choose a vertex $\bar f_0\in V(\cX)$ corresponding to a chamber in $X$. Denote the vertices corresponding to vertices of the chamber by $\bar v_1$, $\bar v_2$ and $\bar v_3$ and the vertices corresponding to panels of the chamber by $\bar e_1$, $\bar e_2$ and $\bar e_3$. Choose $q$ chambers representing the other $\Gamma$-orbits of chambers and denote the corresponding vertices in $\cX$ by $\bar f_1,\ldots,\bar f_q$.

    Now we construct a quotient complex of groups following \cite[III.$\cC$.2.9]{BH:NPC:99}. For this, note that, since the $\Gamma$-action is regular on panels, the vertices we just named form a system of representatives for all orbits of vertices of $\cX$. It is not hard to see that the quotient scwol then has the following structure
    \begin{align*}
       V(\Gamma\backslash\!\backslash \cX) &\coloneq \{ v_1, v_2, v_3\} \sqcup \{e_1,e_2,e_3\} \sqcup \{ f_j : j \in J \} \\
       E(\Gamma\backslash\!\backslash \cX) &\coloneq \{ v_\alpha \leftarrow e_\beta: \alpha\neq\beta\} \sqcup \{ v_\alpha \leftarrow f_j : j\in J\} \sqcup \{e_\beta \leftarrow f_j : j\in J\},
    \end{align*}
where we denote the orbits of vertices by removing the bar. The quotient scwol is also shown in Figure \ref{fig:quotient_scwol} for $q=2$.

    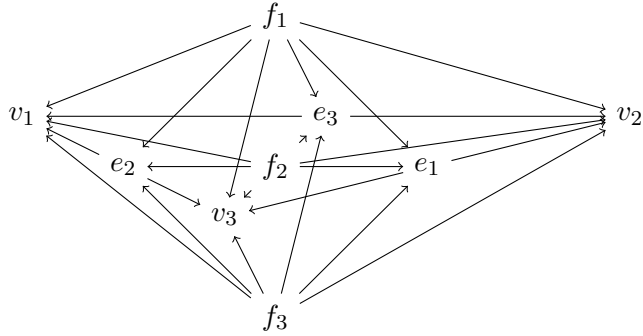
\begin{figure}[hbt]
    \centering
\begin{tikzpicture}[scale=4,font=\small]
	\node (p1) at (0,0,0) {$v_1$};
	\node (p2) at (2,0,0) {$v_2$};
	\node (p3) at (1,0,.866) {$v_3$};
	\node (l3) at (1,0,0) {$e_3$} edge[->] (p1) edge[->] (p2);
	\node (l1) at (1.5,0,.433) {$e_1$} edge[->] (p2) edge[->] (p3);
	\node (l2) at (.5,0,.433) {$e_2$} edge[->] (p1) edge[->] (p3);
	\node (f1) at (1,.5,.433) {$f_1$} edge[->] (p1) edge[->] (p2) edge[->] (p3) edge[->] (l1) edge[->] (l2) edge[->] (l3);
	\node (f2) at (1,0,.433) {$f_2$} edge[->] (p1) edge[->] (p2) edge[->] (p3) edge[->] (l1) edge[->] (l2) edge[->] (l3);
	\node (f3) at (1,-.5,.433) {$f_3$} edge[->] (p1) edge[->] (p2) edge[->] (p3) edge[->] (l1) edge[->] (l2) edge[->] (l3);
\end{tikzpicture} \caption{The quotient scwol $\Gamma\backslash\!\backslash\cX$ for $q=2$}\label{fig:quotient_scwol}
\end{figure}

The vertex groups are just the stabilisers. All of these are hence trivial except the stabilisers $\Gamma_{\bar v_\alpha}$ which are the Singer groups $S_\alpha$. All monomorphisms are obviously trivial.

It remains to determine the twist elements as in \cite[III.$\cC$.2.9]{BH:NPC:99}. For every edge $a\in E(\Gamma\backslash\!\backslash \cX)$, there is precisely one preimage $\bar a\in E(\cX)$ satisfying $i(\bar a) = \overline{i(a)}$, but in general $t(\bar a)\neq \overline{t(a)}$. We choose $h_a$ such that $h_a(t(\bar a)) = \overline{t(a)}$.

We first consider the edges of the form $v_\alpha\leftarrow e_\beta$. For these, we obviously have $h_{v_\alpha\leftarrow e_\beta}=1$. Now consider the edges of the form $e_\beta\leftarrow f_j$. There is exactly one preimage $\bar e_\beta'\leftarrow \bar f_j$, where $\bar e_\beta'$ is in the orbit $e_\beta$. Since the $\Gamma$-action is regular on panels, there is exactly one element $h_{e_\beta\leftarrow f_j}\in\Gamma$ satisfying
\[
h_{e_\beta\leftarrow f_j}(t(\bar e_\beta'\leftarrow \bar f_j) ) = \bar e_\beta.
\]

For edges of the type $\bar v_\alpha' \leftarrow \bar f_j$, where $\bar v_\alpha'$ is in the orbit $v_\alpha$, set $h_{v_\alpha\leftarrow f_j}=h_{e_\beta\leftarrow f_j}$ where $\beta \equiv_3 \alpha + 1$. Since by construction $g_{a,b}= h_ah_b h_{ab}^{-1}$, we have
\[
g_{v_\alpha\leftarrow e_\beta\leftarrow f_j} = \underbrace{h_{v_\alpha\leftarrow e_\beta}}_{=1} h_{e_\beta\leftarrow f_j}h_{v_\alpha\leftarrow f_j}^{-1}.
\]
This implies
\[
g_{v_\alpha\leftarrow e_\beta \leftarrow f_j}= 1\quad\text{for all } \beta-\alpha\equiv_3 1,\, j\in J.
\]
The other twist elements then necessarily form three ordered difference sets $d_\alpha$ as follows
\[
g_{v_\alpha\leftarrow e_\beta \leftarrow f_j}= d_\alpha(j) \quad\text{for } \beta-\alpha\equiv_3 2,\, j\in J,
\]
since the element $h_{e_\beta\leftarrow f_j}h^{-1}_{e_{\beta-1}\leftarrow f_j}$ transports the unique vertex in the orbit $f_j$ adjacent to $\bar e_{\beta-1}$ to the unique vertex in the orbit $f_j$ adjacent to $\bar e_\beta$. Since the vertices $\{\bar v_1,\bar v_2,\bar v_3,\bar e_1,\bar e_2,\bar e_3,\bar f_0\}$ span a subcomplex, we have
\[
d_1(0)=d_2(0)=d_3(0)=1.
\]
By the same calculation as in the proof of Proposition \ref{prop:fundamental_group_a2} we obtain the required presentation.
\end{Proof}

\begin{Remark}
    In particular, this means that any lattice constructed as in section \ref{subsec:general_a2_construction} admits a presentation as in Theorem \ref{th:a2_classification}. This will simplify our considerations in the following sections.
\end{Remark}

\begin{Corollary}\label{cor:a2_building_description}
    In this situation, the building $X$ is isomorphic to the flag complex over the graph with vertices
    \[
        V(X) = \Gamma/S_1 \sqcup \Gamma/S_2 \sqcup \Gamma/S_3
    \]
    and edges
    \[
        E(X) = \{ (gS_1,gS_2), (gS_2,gS_3), (gS_3,gS_1) : g\in \Gamma\}
    \]
\end{Corollary}

\subsection{Lattices generated by cyclic Singer groups}\label{subsec:cyclic_lattices}

In this section, we start with cyclic Singer groups to obtain very simple lattices. Let $\cI$ be the classical projective plane of order $q$. Take three cyclic Singer groups $S_\alpha$ of $\cI$ of order $q^2+q+1$ with respective generators $\sigma_\alpha$. We choose the points and lines for the construction of the difference sets $D_\alpha$ to be incident, such that $1\in D_\alpha$. We write
\[
	\Delta_\alpha \coloneq \{ \delta\in\Z/(q^2+q+1) : \sigma_\alpha^\delta \in D_\alpha \}.
\]
The sets of numbers $\Delta_\alpha$ are then \emph{difference sets} in the classical sense. In addition, we choose bijections $\delta_\alpha:J = \{0,1,\ldots,q\}\rightarrow \Delta_\alpha$ satisfying $\delta_\alpha(0)=0$. If we apply the construction from section \ref{subsec:general_a2_construction}, we obtain the following very simple presentation.

\begin{Theorem}[Lattices generated by cyclic Singer groups]\label{th:cyclic_lattices}
	For any prime power $q$ and any three classical difference sets $\Delta_1$, $\Delta_2$ and $\Delta_3$ containing 0, and for any bijections $\delta_\alpha:\{0,1,\ldots,q\}\rightarrow \Delta_\alpha$ satisfying $\delta_\alpha(0)=0$, the group $\Gamma$ with presentation
	\[
		\Gamma = \bigl\langle \sigma_1, \sigma_2, \sigma_3 \,\big|\, \sigma_1^{q^2+q+1}=\sigma_2^{q^2+q+1}=\sigma_3^{q^2+q+1}=1,\, \sigma_1^{\delta_1(j)}\sigma_2^{\delta_2(j)}\sigma_3^{\delta_3(j)} = 1 \quad\forall j \in J\bigr\rangle
\]
	is a uniform lattice in a building of type $\tilde A_2$.
\end{Theorem}

\begin{Examples}
	We give explicit presentations for three lattices in the two smallest cases. More difference sets $\Delta$ can be obtained from \cite{LaJolla}.
	\begin{enumerate}
		\item For $q=2$, we choose all difference sets to be $\Delta=\{0,1,3\}$. We obtain the lattice
			\[
				\Gamma_2 \coloneq \langle \sigma_1,\sigma_2,\sigma_3 \,|\, \sigma_1^7=\sigma_2^7=\sigma_3^7 = \sigma_1\sigma_2\sigma_3 = \sigma_1^3\sigma_2^3\sigma_3^3 = 1\rangle.
			\]
			It can easily be seen that $H_1(\Gamma_2)\cong (\Z/7)^2$.
		\item By changing the order of the first difference set, we obtain a new lattice which is not isomorphic to the first one.
			\[
				\Gamma_2' \coloneq \langle \sigma_1,\sigma_2,\sigma_3 \,|\, \sigma_1^7=\sigma_2^7=\sigma_3^7 = \sigma_1^3\sigma_2\sigma_3 = \sigma_1\sigma_2^3\sigma_3^3 = 1\rangle.
			\]
			For the abelianisation, we obtain $H_1(\Gamma_2') \cong \Z/7$.
		\item For $q=3$, we choose the difference set $\Delta=\{0,1,3,9\}$. We obtain
			\[
			\Gamma_3 \coloneq \langle \sigma_1,\sigma_2,\sigma_3 \,|\, \sigma_1^{13}=\sigma_2^{13}=\sigma_3^{13} = \sigma_1\sigma_2\sigma_3 = \sigma_1^3\sigma_2^3\sigma_3^3 =\sigma_1^9\sigma_2^9\sigma_3^9 = 1\rangle.
			\]
			Here, we have $H_1(\Gamma_3)\cong (\Z/13)^2$.
	\end{enumerate}
    For $q=5$, it is possible to construct perfect lattices with this method.
\end{Examples}

\begin{Remark}
	In \cite{KMW:A2l:84}, Köhler, Meixner and Wester showed that the following group is a  chamber-regular lattice in the building associated to $\Sl_3(\F_2(\!(t)\!))$:
    \[
    \Gamma = \langle a,b,c\,|\, a^3=b^3=c^3=1,\, (ab)^2=ba, (ac)^2=ca, (c^{-1}b)^2=bc^{-1} \rangle.
    \]
    It is not hard to verify that the subgroup generated by $a^{-1}b^{-1}$, $bc^{-1}$ and $ca$ is isomorphic to the lattice $\Gamma_2$ from the example above.

    Computer searches done by the author using \cite{SAGE} and \cite{GAP4} have not yielded any embeddings of our lattices in $\Sl_3(\F_q(\!(t)\!))$ for $q>2$ yet.
\end{Remark}

\subsection{Spheres of radius two}\label{subsec:spheres}

For any vertex of an affine building of type $\tilde A_2$, the combinatorial sphere of radius one, usually called the link, forms a projective plane. The spheres of larger radii also admit interesting incidence structures, the so-called \emph{Hjelmslev planes}, see \cite{HvM:PHP:89}. It is a well-known fact that the two families of classical buildings of type $\tilde A_2$ associated to the groups $\PSl_3(\Q_p)$ and $\PSl_3(\F_p(\!(t)\!))$ can already be distinguished by inspecting spheres of radius two, respectively Hjelmslev planes of level two.

We will use a criterion by Cartwright, Mantero, Steger and Zappa, see \cite[section 8]{CMSZ2}, to distinguish between these two affine buildings. Consider the following construction.

Denote the type set of a building $X$ of type $\tilde A_2$ by \{1,2,3\}. Fix a chamber $c_0$ and denote its vertices of type $1$, $2$ and $3$ by $v_1$, $v_2$ and $v_3$, respectively. Denote by $\cP$ and $\caL$ the sets of vertices adjacent to $v_1$ of types $3$ and $2$, respectively. Then of course $\cP$ and $\caL$ with the adjacency relation in the building form a projective plane.

Now consider the vertices of type $2$ and $3$ at distance $2$ from the vertex $v_1$ and denote these sets by $\cP^2$ and $\caL^2$, respectively.

\begin{figure}[hbt]
\centering
\begin{tikzpicture}[scale=0.8]
    \tikzstyle{every node}=[circle, draw, inner sep=0pt, minimum width=8pt]
    \draw (0,0) node[fill=black] {} -- (4,0) node[fill=white] {} -- (8,0) node[fill=gray] {};
    \draw (0,0) -- (2,3.464) node[fill=gray] {} -- (4,6.928) node[fill=white] {};
    \draw (2,3.464) node[fill=gray] {} -- (4,0) node[fill=white] {} -- (6,3.464) -- cycle;
    \draw (4,6.928) node[fill=white] {} -- (6,3.464) node[fill=black] {} -- (8,0) node[fill=gray] {};
    \tikzstyle{every node}=[]
    \node (v1) at (-.7,0) {$v_1$};
    \node (v2) at (8.75,0) {$x_2$};
    \node (v3) at (4.75,6.928) {$x_3$};
\end{tikzpicture} \caption{Adjacency in the Hjelmslev plane of level 2}\label{fig:adjacency}
\end{figure}
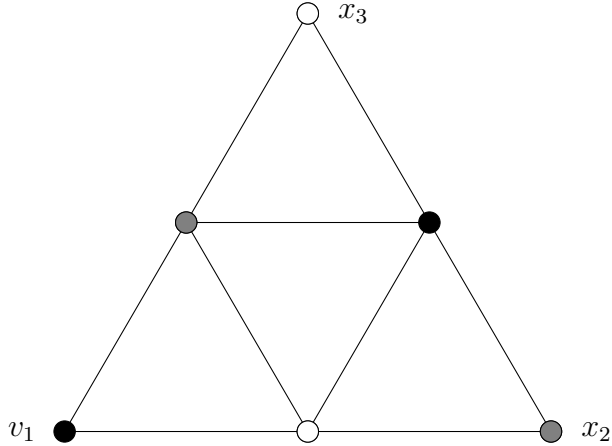

\begin{Definition}
    Two vertices $x_2 \in \cP^2$, $x_3\in \caL^2$ are \emph{adjacent}, $x_2\sim x_3$, if the configuration in Figure \ref{fig:adjacency},
    where the colours of the vertices indicate the type, exists in the building. The incidence structure $(\cP^2,\caL^2,\sim)$ is called a \emph{Hjelmslev plane} of level two.
\end{Definition}

\noindent We will give an explicit construction of the Hjelmslev plane in the case of a building $X$ of type $\tilde A_2$ admitting a panel-regular lattice $\Gamma$.

Denote the vertex stabilisers of the vertices $v_\alpha$ by $S_\alpha$, these are Singer groups in the corresponding vertex links. For every vertex $v_\alpha$, the link $\lk_X(v_\alpha)$ is a projective plane of order $q$, where call the vertices of type $\alpha+1 \pmod{3}$ lines and the vertices of type $\alpha+2 \pmod{3}$ points.

 For each of these projective planes, we construct a difference set $D_\alpha$ with respect to the points and lines given by the other vertices of $c_0$, respectively. We take the description of the building from Corollary \ref{cor:a2_building_description} and prove a simple lemma.

\begin{Lemma}\label{la:triangles}
    We have
    \begin{itemize}
        \item If $h,f\in S_1$ and $(v_1,hv_2,fv_3)$ is a triangle in $X$, then $f^{-1}h\in D_1$.
        \item If $g,f\in S_2$ and $(gv_1,v_2,fv_3)$ is a triangle in $X$, then $g^{-1}f\in D_2$.
        \item If $g,h\in S_3$ and $(gv_1,hv_2,v_3)$ is a triangle in $X$, then $h^{-1}g\in D_3$.
    \end{itemize}
\end{Lemma}

\begin{Proof}
    We will only prove the first claim, the other two are analogous. If $(v_1, hv_2,fv_3)$ is a triangle, then so is $(v_1,f^{-1}hv_2,v_3)$ which we obtain by multiplying with $f^{-1}$. Since $S_1$-translates of $v_2$ correspond to lines in the projective plane $\lk_X(v_1)$ by the definition of the difference set $D_1$, the line $f^{-1}h v_2$ is incident to the point $v_3$ if and only if $f^{-1}h\in D_1$.
\end{Proof}

\begin{Lemma}\label{lem:points_and_lines}
    We have
    \[
    \cP = \{ s_1 v_3 : s_1 \in S_1 \},\qquad \caL = \{ t_1 v_2 : t_1 \in S_1 \},
    \]
    as well as
    \[
    \cP^2 = \{ s_1s_3 v_2: s_1\in S_1, s_3\in S_3\backslash D_3^{-1} \},\qquad \caL^2 = \{ t_1t_2 v_3 : t_1\in S_1, t_2\in S_2\backslash D_2 \}.
    \]
\end{Lemma}

\begin{Proof}
    The first claim is obvious, since $S_1$ acts regularly on points and lines of the projective plane $\lk(v_1)$. The second claim is a simple consequence of Lemma \ref{la:triangles} --- the vertex $t_2v_3$ is not adjacent to $v_1$ if and only if $t_2\not\in D_2$, and similarly for $s_3$.
\end{Proof}

\noindent Adjacency in the Hjelmslev plane is relatively complex to describe but becomes manageable for cyclic Singer groups. The following observation will be used frequently in this section:

From the presentation of $\Gamma$ in Theorem \ref{th:a2_classification}, we see that for $d_1\in D_1$, there are always elements $d_2\in D_2$, $d_3\in D_3$ such that $d_1d_2d_3=1$.

\begin{Lemma}\label{lem:hjelmslev_adjacency}
    Fix a point $s_1s_3v_2\in\cP^2$ and a line $t_1t_2v_3\in\caL^2$, where $s_1,t_1\in D_1$, $t_2\in D_2$ and $s_3\in D_3$. In the Hjelmslev plane of level 2 around $v_1$, we have $s_1s_3v_2 \sim t_1t_2v_3$ if and only if the following conditions hold.
    \begin{description}
        \item[(C1)] We have $s_1^{-1}t_1 \in D_1$.
    \end{description}
     There are hence elements $d_2\in D_2$, $d_3\in D_3$ such that $(s_1^{-1}t_1)d_2d_3=1$.
     \begin{description}
        \item[(C2)] There is an element $n_2\in t_2D_2^{-1} \cap d_2D_2^{-1}$ such that
            \[
            s_3^{-1}d_3^{-1}e_3 \in D_3,
            \]
            where $e_1(n_2^{-1}d_2)e_3 = 1$ for some elements $e_1\in D_1$ and $e_3\in D_3$.
    \end{description}
    Note that the elements $e_1$ and $e_3$ exist since $n_2\in d_2D_2^{-1}$.

    If the lattice arises from cyclic Singer groups $S_\alpha=\langle \sigma_\alpha\rangle$, $\alpha\in\{1,2,3\}$, as in section \ref{subsec:cyclic_lattices} and if the associated unordered difference sets $\Delta_\alpha=\Delta$ are all equal, we obtain:
    \[
    \sigma_1^{j_1}\sigma_3^{j_3}v_2 \sim \sigma_1^{k_1}\sigma_2^{k_2}v_3\quad\Leftrightarrow\quad\begin{cases}
        k_1 - j_1 \in\Delta \text{ and}\\
        \exists n\in (k_2 - \Delta) \cap (-j_3-\Delta) \cap ( k_1-j_1-\Delta).
    \end{cases}
    \]
\end{Lemma}

\begin{Proof}
	Assume that $s_1s_3v_2\sim t_1t_2v_3$. Then there must be an element $n_2\in S_2$ such that we have the configuration of Figure \ref{fig:proof_adjacency} in $X$. We will now investigate the triangles in the order given by roman numerals and apply Lemma \ref{la:triangles} repeatedly to obtain the required relations.
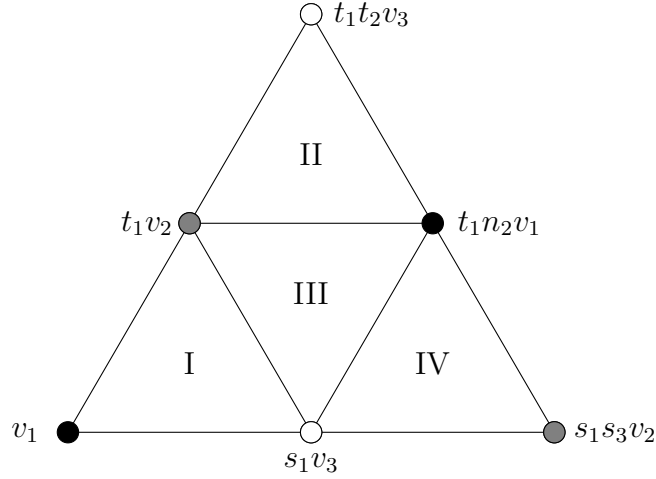
\begin{figure}[hbt]
    \centering
    \begin{tikzpicture}[scale=0.8]
        \tikzstyle{every node}=[circle, draw, inner sep=0pt, minimum width=8pt]
        \draw (0,0) node[fill=black] {} -- (4,0) node[fill=white] {} -- (8,0) node[fill=gray] {};
        \draw (0,0) -- (2,3.464) node[fill=gray] {} -- (4,6.928) node[fill=white] {};
        \draw (2,3.464) node[fill=gray] {} -- (4,0) node[fill=white] {} -- (6,3.464) -- cycle;
        \draw (4,6.928) node[fill=white] {} -- (6,3.464) node[fill=black] {} -- (8,0) node[fill=gray] {};
        \tikzstyle{every node}=[]
        \node (v1) at (-.7,0) {$v_1$};
        \node (v2) at (9,0) {$s_1s_3v_2$};
        \node (v3) at (5,6.928) {$t_1t_2v_3$};
        \node (v4) at (4,-.5) {$s_1v_3$};
        \node (v5) at (1.3,3.464) {$t_1v_2$};
        \node (v6) at (7.1,3.464) {$t_1n_2v_1$};
        \node (i) at (2,1.154) {I};
        \node (ii) at (4,4.618) {II};
        \node (iii) at (4,2.309) {III};
        \node (iv) at (6,1.154) {IV};
    \end{tikzpicture} \caption{The configuration in the proof of Lemma \ref{lem:hjelmslev_adjacency}}\label{fig:proof_adjacency}
    \end{figure}
    Since $(v_1,t_1v_2,s_1v_3)$ form triangle I, we obtain
    \begin{equation}\label{eq:I}
        s_1^{-1}t_1\in D_1.
    \end{equation}
    There are hence elements $d_2\in D_2$ and $d_3\in D_3$ such that $(s_1^{-1}t_1)d_2d_3=1$. By looking at triangle II translated by $t_1^{-1}$, we obtain
    \begin{equation}\label{eq:II}
        n_2^{-1}t_2\in D_2.
    \end{equation}

    \noindent From triangle III, by observing that $s_1v_3=s_1d_3^{-1}v_3$ and by translating by $t_1^{-1}$, we obtain the triangle \[
    (n_2v_1,v_2,t_1^{-1}s_1d_3^{-1}v_3)=(n_2v_1,v_2,d_2v_3),
    \]
    which again by Lemma \ref{la:triangles} yields
    \begin{equation}
        n_2^{-1}d_2 \in D_2,\label{eq:III}
    \end{equation}
    and there are hence elements $e_1\in D_1$ and $e_3\in D_3$ such that $e_1(n_2^{-1}d_2)e_3=1$. For the last triangle IV, we first write $t_1n_2v_1=t_1n_2e_1^{-1}v_1$. By translating by $s_1^{-1}$, we obtain the triangle
    \begin{align}
        (s_1^{-1}t_1n_2e_1^{-1}v_1,s_3v_2,v_3)\nonumber \\
        \intertext{Substituting $s_1^{-1}t_1$ by $d_3^{-1}d_2^{-1}$ yields the triangle}
        (d_3^{-1}d_2^{-1}n_2e_1^{-1}v_1,s_3v_2,v_3)\nonumber \\
        \intertext{Finally, we substitute $d_2^{-1}n_2=e_3e_1$ and obtain the triangle}
        (d_3^{-1}e_3v_1,s_3v_2,v_3).\nonumber \\
        \intertext{We can now apply Lemma \ref{la:triangles} again to obtain}
        s_3^{-1}d_3^{-1}e_3 \in D_3. \label{eq:IV}
    \end{align}
    By assembling equations \eqref{eq:I},\eqref{eq:II}, \eqref{eq:III} and \eqref{eq:IV}, we obtain conditions (C1) and (C2). For the reverse direction, observe that (C1) and (C2) imply directly that the configuration of Figure \ref{fig:proof_adjacency} exists in the building.

    The second claim is then a simple calculation.
\end{Proof}

\noindent There is an obvious projection $\psi:(\cP^2,\caL^2,\sim) \rightarrow (\cP,\caL,\cF)$ given by

\begin{align*}
    \psi: \cP^2 &\rightarrow \cP: s_1s_3v_2 \mapsto s_1v_3\\
         \caL^2 &\rightarrow \caL: t_1t_2v_3 \mapsto t_1v_2.
\end{align*}

\begin{Lemma}[Lemma 8.1 in \cite{CMSZ2}]\label{l:cmsz}
    For $p,p'\in \cP^2$ consider the images $\psi(p)$ and $\psi(p')$. If $\psi(p)\neq \psi(p')$, there is a unique $l\in \caL^2$ such that $p\sim l\sim p'$. If $\psi(p)=\psi(p')$, there are $q$ distinct such lines. The same is true if the roles of points and lines are reversed.
\end{Lemma}

\begin{Proposition}\label{prop:hjelmslev_splitting}
    If the lattice $\Gamma$ and the building $X$ arise from cyclic Singer groups as in section \ref{subsec:cyclic_lattices} and if in addition all difference sets $\Delta_\alpha=\Delta$ are equal, there is a splitting map $\iota: (\cP,\caL,\cF) \rightarrow (\cP^2,\caL^2,\sim)$ satisfying $\psi\circ\iota=\id_{(\cP,\caL,\cF)}$.
\end{Proposition}

\begin{Proof}
    With the same notation as in Lemma \ref{lem:hjelmslev_adjacency}, we choose an element
    \[
    m\in(\Z/(q^2+q+1))\setminus(\Delta\cup(-\Delta)).
    \]
    This is always possible for cardinality reasons, since $q^2+q+1 > 2q+2$ for $q\geq 2$. The map $\iota$ is given by
    \begin{align*}
        \cP & \rightarrow \cP^2&&\text{and}             & \caL & \rightarrow \caL^2 \\
        s_1v_3 &\mapsto s_1\sigma_3^{-m}v_2 &&& t_1v_2 & \mapsto t_1\sigma_2^{m}v_3.
    \end{align*}
    This map preserves incidence by Lemma \ref{lem:hjelmslev_adjacency}, since for any $j_1, k_1$, the set
    \[
        (m - \Delta) \cap ( k_1 - j_1 - \Delta)
    \]
    has exactly one element, which corresponds to the intersection of lines in the projective plane $\lk_X(v_1)$.
\end{Proof}

\begin{Definition}[Section 8 in \cite{CMSZ2}]
    Consider $\cP'\subset \cP^2$ and $\caL'\subset \caL^2$. We say that $(\cP',\caL')$ is a \emph{substructure} of $(\cP^2,\caL^2)$ if for all $p,p'\in\cP'$ satisfying $\psi(p)\neq\psi(p')$, the unique line $l\in\caL^2$ incident with both $p$ and $p'$ is already contained in $\caL'$ and if the same condition holds with points and lines exchanged. The substructure \emph{generated by a set $\cP''\subset \cP^2$} is the smallest substructure containing $\cP''$.
\end{Definition}

\noindent We will use the following characterisation to distinguish spheres of radius $2$ in buildings of type $\tilde A_2$ by Cartwright-Mantero-Steger-Zappa in \cite{CMSZ2}.

\begin{Proposition}[Proposition 8.6 in \cite{CMSZ2}]\label{prop:cmsz}
   Assume that $X$ is the building associated to the projective special linear group $\PSl_3(K)$, where $K$ is a local, non-Archimedean field with residue field of prime order $p$. Choose a vertex $v_1\in X$ and construct the projective plane $(\cP,\caL,\cF)$ and the Hjelmslev plane of level 2 denoted by $(\cP^2,\caL^2,\sim)$ as above. Consider four points $p_1,\ldots,p_4\in\cP^2$ such that no three points of $\psi(p_1),\ldots,\psi(p_4)$ are collinear.
   \begin{itemize}
       \item If $K=\Q_p$, the substructure of the Hjelmslev plane generated by the set of points $\{p_1,\ldots,p_4\}$ is all of $(\cP^2,\caL^2,\sim)$.
       \item Otherwise, the substructure generated by the set of points $\{p_1,\ldots,p_4\}$ is a projective plane of order $p$.
   \end{itemize}
\end{Proposition}

\begin{Corollary}
    The panel-regular lattices constructed out of cyclic Singer groups with identical difference sets cannot be contained in the building associated to $\PSl_3(\Q_p)$.
\end{Corollary}

\begin{Proof}
    Take any four points $p_1,\ldots,p_4$ in $(\cP,\caL,\cF)$ such that no three points are collinear. Then $\iota(p_1),\ldots, \iota(p_4)$ satisfy the conditions of Proposition \ref{prop:cmsz}. But by Lemma \ref{l:cmsz} and by Proposition \ref{prop:hjelmslev_splitting}, the substructure generated by these four points is exactly $\im(\iota)$, which is a projective plane of order $p$. By Proposition \ref{prop:cmsz}, if the building $X$ is associated to $\PSl_3(K)$, then necessarily $K\neq \Q_p$.
\end{Proof}

\begin{Remark}
    Conjecturally, these very simple lattices should be contained in $\Sl_3(\F_p(\!(t)\!))$. It might be possible to prove this by investigating the structure of all Hjelmslev planes of all levels and prove a splitting lemma for every step. This seems tedious and difficult, however. Surprisingly, we have not found any realisation of such a lattice in $\Sl_3(\F_p(\!(t)\!))$ for $p\neq 2$ yet.
\end{Remark}

\section{Lattices in buildings of type \texorpdfstring{$\tilde C_2$}{\textasciitilde C2}}\label{sec:c2}

In this section, we will construct lattices in buildings of type $\tilde C_2$. This will be done using the slanted symplectic quadrangles discussed in section \ref{subsec:GQs}. We will give two different constructions of panel-regular lattices in buildings of type $\tilde C_2$.

\paragraph{Situation} Fix a prime power $q>2$, and set $J=\{0,1,\ldots, q+1\}$. Let $\cI$ and $\cI'$ be two copies of the slanted symplectic quadrangle of order $(q-1,q+1)$. For the quadrangle $\cI$, we fix a point $p$ and denote the set of lines through $p$ by $L$. The set of flags $F$ is given by $F=\{ (p,l):l\in L\}$. The associated Singer group will be denoted by $S$, the line stabilisers by $S_l$ for lines $l\in L$. We fix the same objects for $\cI'$ and add a prime $'$ to the notation.

\subsection{Lattices acting regularly on two types of panels}\label{subsec:c2_two_panels}

Since $q$ determines the quadrangle and the Singer group uniquely up to isomorphism, there is a bijection between the sets of line representatives $L$ and $L'$ and the corresponding line stabilisers $S_l$ and $S'_{l'}$ are pairwise isomorphic.

\paragraph{Construction} We choose two bijections $\lambda:J\rightarrow L$, $\lambda':J\rightarrow L'$ and an abstract group $S_j\cong \Z/q$ with isomorphisms $\psi_j: S_j\rightarrow S_{\lambda(j)}$ and $\psi'_j: S_j\rightarrow S'_{\lambda'(j)}$ for every $j\in J$. Consider the complex of groups $G(\cY)$ over the scwol $\cY$ with vertices
\begin{align*}
	V(\cY) \coloneq&\, \{v,v',w\} \sqcup \{e,e'\} \sqcup \{e_j : j\in J\} \sqcup \{f_j : j\in J \} \\
	\intertext{and edges}
	E(\cY) \coloneq&\, \{w \leftarrow e, w\leftarrow e', v\leftarrow e, v'\leftarrow e'\} \sqcup \{ v\leftarrow e_j, v'\leftarrow e_j : j\in J\} \\
	& \quad\sqcup \{ v\leftarrow f_j, v'\leftarrow f_j, w\leftarrow f_j : j\in J \} \sqcup \{e\leftarrow f_j, e'\leftarrow f_j, e_j\leftarrow f_j : j \in J\}.
\end{align*}
Figure \ref{fig:cYfor3} illustrates the scwol $\cY$ in the case $q=3$.
\begin{figure}[hbt]
\centering
\begin{tikzpicture}[scale=4,font=\small]
	\node (p1) at (0,0,0) {$w$};
	\node (p2) at (2,0,0) {$v'$};
	\node (p3) at (1,0,.866) {$v$};
	\node (l3) at (1,0,0) {$e'$} edge[->] (p1) edge[->] (p2);
	\node (l2) at (.5,0,.433) {$e$} edge[->] (p1) edge[->] (p3);
	\node[color=gray] (la) at (1.5,1,.433) {$e_0$};
	\node[color=gray] (lb) at (1.5,.5,.433) {$e_1$};
	\node (lc) at (1.5,0,.433) {$e_2$} edge[->] (p2) edge[->] (p3);
	\node[color=gray] (ld) at (1.5,-.5,.433) {$e_3$};
	\node[color=gray] (le) at (1.5,-1,.433) {$e_4$};
	\node[color=gray] (fa) at (1,.66,.433) {$f_0$} edge[->,color=lightgray] (p1) edge[->,color=lightgray] (p2) edge[->,color=lightgray] (p3) edge[->,color=lightgray] (la) edge[->,color=lightgray] (l2) edge[->,color=lightgray] (l3);
	\node[color=gray] (fb) at (1,.33,.433) {$f_1$} edge[->,color=gray] (p1) edge[->,color=gray] (p2) edge[->,color=gray] (p3) edge[->,color=gray] (lb) edge[->,color=gray] (l2) edge[->,color=gray] (l3);
	\node (fc) at (1,0,.433) {$f_2$} edge[->] (p1) edge[->] (p2) edge[->] (p3) edge[->] (lc) edge[->] (l2) edge[->] (l3);
	\node[color=gray] (fd) at (1,-.33,.433) {$f_3$} edge[->,color=gray] (p1) edge[->,color=gray] (p2) edge[->,color=gray] (p3) edge[->,color=gray] (ld) edge[->,color=gray] (l2) edge[->,color=gray] (l3);
	\node[color=gray] (fe) at (1,-0.66,.433) {$f_4$} edge[->,color=lightgray] (p1) edge[->,color=lightgray] (p2) edge[->,color=lightgray] (p3) edge[->,color=lightgray] (le) edge[->,color=lightgray] (l2) edge[->,color=lightgray] (l3);
	\draw[->,color=lightgray] (la) .. controls (1.75,.7,.216) .. (p2);
	\draw[->,color=lightgray] (la) .. controls (1.25,.7,.649) .. (p3);
	\draw[->,color=lightgray] (le) .. controls (1.75,-.7,.216) .. (p2);
	\draw[->,color=lightgray] (le) .. controls (1.25,-.7,.649) .. (p3);
	\draw[->,color=gray] (lb) .. controls (1.75,.35,.216) .. (p2);
	\draw[->,color=gray] (lb) .. controls (1.25,.35,.649) .. (p3);
	\draw[->,color=gray] (ld) .. controls (1.75,-.35,.216) .. (p2);
	\draw[->,color=gray] (ld) .. controls (1.25,-.35,.649) .. (p3);
\end{tikzpicture} \caption{The scwol $\cY$ for $q=3$}\label{fig:cYfor3}
\end{figure}

Now choose the vertex groups to be $G_{v}= S$, $G_{v'}= S'$ and $G_{w}= \langle c\,|\, c^{q+2}=1\rangle$. In addition, set $G_{e_j}= S_j$ for all $j\in J$. All other vertex groups are trivial. The only non-trivial monomorphisms are chosen to be the maps $\psi_j$ and $\psi'_j$ for $j\in J$. All twist elements are trivial except
\[
g_{w\leftarrow e\leftarrow f_j} \coloneq c^j.
\]

\noindent We endow $|\cY|$ with a locally Euclidean metric as follows:

Let $\Delta$ be the geometric realisation of one triangle in the affine Coxeter complex of type $\tilde C_2$. For each $j\in J$, we map the subcomplex spanned by $\{v,v',w,e,e',e_j,f_j\}$ onto the barycentric subdivision of $\Delta$ in the obvious way and pull back the metric. We obtain a locally Euclidean metric on $|\cY|$. In particular, the angles at the vertices $v$ and $v'$ are $\pi/4$, the angle at $w$ is $\pi/2$.

\begin{Proposition}\label{prop:c2_developable}
	The complex of groups $G(\cY)$ is developable.
\end{Proposition}

\begin{Proof}
    The proof is analogous to the one of Proposition \ref{prop:a2_is_developable}. By construction of the complex, the geometric links in the local developments at $v$, $v'$ and $w$ are $\cZ(\cI)$, $\cZ(\cI')$ and the scwol associated to a complete bipartite graph of order $(q+2,q+2)$. In particular, they are CAT(1) by Proposition \ref{prop:key}.

	The local developments of the edges $e$ and $e'$ are the same as the closed stars in $\cY$: $(q+2)$ 2-simplices glued along one edge. For the edges $e_j$, the local development is isometric to $q$ triangles glued along one edge. The local developments of the vertices $f_j$ are just flat triangles.

    By Proposition \ref{prop:geometric_link_nonpos_curved}, the complex of groups $G(\cY)$ is non-positively curved and hence developable by Theorem \ref{th:developability_of_nonpositively_curved_complexes}.
\end{Proof}

\begin{Proposition}\label{prop:c2_presentation}
	The fundamental group $\Gamma=\pi_1(G(\cY))$ admits the following presentation:
	\[
        \Gamma = \bigl\langle S,S',c \,\big|\, \text{ all relations in } S,S',\, c^{q+2}=1,\, c^j\psi_j(s)c^{-j}= \psi'_j(s)\quad \forall j \in J, s\in S_j \bigr\rangle.
	\]
\end{Proposition}

\begin{Proof}
	Consider the following maximal spanning subtree $T$:
	\[
		E(T)=\{w\leftarrow e, w\leftarrow e', v\leftarrow e, v'\leftarrow e'\} \sqcup \{w\leftarrow f_j, e_j\leftarrow f_j: j\in J\}
	\]
	As in the proof of Proposition \ref{prop:fundamental_group_a2}, the group is generated by $S$, $S'$, $c$ and group elements for each edge not contained in $T$. Consider, however Figure \ref{fig:c2_triangle} showing one triangle in $\cY$.
	\begin{figure}
	\centering
	\begin{tikzpicture}[scale=6,font=\small]
		\node (v3) at (0,0) {$w$};
		\node (v2) at (0,1) {$v$};
		\node (v1) at (1,0) {$v'$};
		\node (e1) at (0,.5) {$e$} edge[->] (v3) edge [->] (v2);
		\node (e2) at (.5,0) {$e'$} edge[->] (v3) edge [->] (v1);
        \node (el) at (.5,.5) {$e_j$} edge[->,dotted] node[color=gray] {$1$} (v1) edge[->,dotted] node[color=gray]{$c^j$} (v2);
		\node (fl) at (.33,.33) {$f_j$} edge[->] (v3) edge[->,dotted] node[color=gray]{$1$} (e2) edge[->,dotted] node[color=gray]{$1$} (v1) edge[->,dotted] node[color=gray]{$c^j$} (e1) edge[->,dotted] node[color=gray]{$c^j$} (v2) edge[->] (el);
		\node (m12) at (.11,.61) {$1$};
		\node (m21) at (.61,.11) {$1$};
		\node (m13) at (.11,.29) {$c^j$};
		\node (m32) at (.29,.11) {$1$};
		\node (ml2) at (.27,.61) {$1$};
		\node (ml1) at (.61,.27) {$1$};
	\end{tikzpicture} \caption{One triangle in $\cY$}\label{fig:c2_triangle}
\end{figure}
	As before, edges in $T$ are drawn black, the other edges are drawn dotted. Group elements $g$ written on an edge $a$ indicate that $k_a=g$. From the presentation in Definition \ref{def:fundamental_group}, one can see that all additional relations are of the form stated in the result.
\end{Proof}

\begin{Theorem}\label{th:c2_building}
        The universal cover $\cX$ is a building of type $\tilde C_2$. The fundamental group $\Gamma=\pi_1(G(\cY))$ with presentation
	\[
        \Gamma = \bigl\langle S,S',c \,\big|\, \text{ all relations in } S,S',\, c^{q+2}=1,\, c^j\psi_j(s)c^{-j}= \psi'_j(s)\quad \forall j \in J, s\in S_j \bigr\rangle
	\]
        is hence a uniform lattice in the automorphism group of $\cX$, which is panel-regular on two types of panels.
\end{Theorem}

\begin{Proof}
    The proof is analogous to the one of Theorem \ref{th:a2_is_building}, using Proposition \ref{prop:key} and Theorem \ref{th:recognition}.
\end{Proof}

\begin{Remark}
    Except for the slanted symplectic quadrangle $W(3)^\Diamond$, all Singer quadrangles are necessarily exceptional. Hence, except for possibly this case, all buildings constructed in this way are exceptional buildings of type $\tilde C_2$.

    It is not hard to see that $\langle c\vert c^{q+2}=1\rangle$ can be replaced by any group $C$ of order $q+2$. In the above presentation, the term $c^j$ has to be replaced by an arbitrary bijection $J\rightarrow C$.
\end{Remark}

\noindent Finally, we also have a very explicit description of this building as the flag complex over the following graph
\begin{align*}
        V(X) &\coloneq \Gamma/S \sqcup \Gamma/S' \sqcup \Gamma/\langle c\rangle \\
        E(X) &\coloneq \{ (gS,gS'), (gS,g\langle c\rangle), (gS',g\langle c \rangle) : g\in \Gamma\},
\end{align*}
	which could be used to study these exotic affine buildings.

	Instead, we will focus on the two cases where we obtain relatively simple presentations of these lattices.

\subsubsection{The prime case}

If $q=p$ is prime, we have very simple presentations of the Heisenberg groups $S$ and $S'$. Denote their generators by $x, y$ and $x', y'$, respectively and set $z=[x,y]$ and $z'=[x',y']$. The sets of line representatives $L$ and $L'$ can be parametrised by $\bP\F_p^2\sqcup\{0\}$ as in Theorem \ref{th:slanted_sympl_quadrangle}. We pick the following generators for the line stabilisers
\[
c_{[a:b]}= x^ay^bz^{-\tfrac{1}{2}ab},\quad c_0 = z,\quad c'_{[a:b]}= x'^ay'^bz'^{-\tfrac{1}{2}ab},\quad c'_0 = z'.
\]

\begin{Theorem}\label{th:c2_prime_result1}
    For any odd prime $p$, for $J=\{0,1,\ldots,p+1\}$ and any two bijections $\lambda, \lambda' : J \rightarrow \bP\F_p^2\sqcup\{0\}$, consider the group $\Gamma$ presented by
    \[
    \Gamma = \Biggl\langle x,y,x',y',c \,\Bigg|\, \begin{array}{c} z= xyx^{-1}y^{-1}, z'=x'y'x'^{-1}y'^{-1}, \\
            x^p=y^p=x'^p=y'^p=c^p=z^p=z'^p=1, \\ xz=zx, yz=zy, x'z'=z'x', y'z' = z'y', \\
            c^j c_{\lambda(j)} c^{-j} = c'_{\lambda'(j)} \qquad\forall j\in J \end{array}\Biggr\rangle.
    \]
    Then $\Gamma$ is a uniform lattice in a building of type $\tilde C_2$, where all vertex links which are quadrangles are isomorphic to the slanted symplectic quadrangle of order $(p-1,p+1)$.
\end{Theorem}

\begin{Proof}
    This is a direct application of Theorem \ref{th:c2_building}.
\end{Proof}

\subsubsection{The abelian case}

If $q$ is even, then the associated Singer groups $S$ and $S'$ are isomorphic to $\F_q^3$. We write $J=\{0,1,\ldots,q+1\}$ and fix two bijections
\[
\lambda, \lambda': J \rightarrow \bP\F_q^2\sqcup\{0\}.
\]
We write $S_{[a:b]}=\F_q(a,b,0)^T$ and $S_{0}= \F_q(0,0,1)^T$ for the line stabilisers in $S$ as in Theorem \ref{th:slanted_sympl_quadrangle} and add a prime $'$ for the respective groups in $S'$. Finally, we fix isomorphisms of abstract groups $\psi_j,\psi_j': \Z/q \rightarrow S_{\lambda(j)},S'_{\lambda'(j)}$ for any $j\in J$.

\begin{Theorem}
    Consider the group $\Gamma$ given by
    \[
    \Gamma= \bigl(S * S' * \langle c \rangle \bigr) / \langle c^{q+2}=1, c^j (\psi_j(x)) c^{-j} = \psi'_j(x) : x\in\Z/q, j\in J\rangle.
    \]
    Then $\Gamma$ is a uniform lattice in an exotic building of type $\tilde C_2$, where all vertex links which are quadrangles are isomorphic to the slanted symplectic quadrangle of order $(q-1,q+1)$.
\end{Theorem}

\begin{Proof}
	Again, this is a simple application of Theorem \ref{th:c2_building}.
\end{Proof}

\subsection{Lattices acting regularly on one type of panel}\label{subsec:c2_one_panel}

In this section, we will consider lattices acting regularly on only one type of panel in a building of type $\tilde C_2$. We will concentrate on the simpler case where the type of this panel corresponds to the two extremal vertices in the $\tilde C_2$-diagram.

As before, let $\cI$ and $\cI'$ be two slanted symplectic quadrangles of the same order $(q-1,q+1)$ and let $S$ and $S'$ be the associated Singer groups. Pick points $p$ and $p'$ and denote the sets of lines incident to $p$ and $p'$ by $L$ and $L'$, respectively. Then $F=\{(p,l):l\in L\}$ and $F'=\{(p',l') : l'\in L'\}$ are sets of flag representatives.

\paragraph{Construction} We set $J=\{0,1,\ldots,q+1\}$ and choose two bijections $\lambda:J\rightarrow L$, $\lambda':J\rightarrow L'$. We consider the scwol $\cY$ given by
\begin{align*}
    V(\cY) &\coloneq \{ v,v'\} \sqcup \{v_j : j\in J\} \sqcup \{e\} \sqcup \{e_j, e'_j : j\in J\} \sqcup \{f_j : j\in J\} \\
    E(\cY) &\coloneq \{ v\leftarrow e, v'\leftarrow e\} \sqcup \{v \leftarrow e_j, v_j \leftarrow e_j: j\in J\} \sqcup \{v'\leftarrow e'_j, v_j\leftarrow e'_j : j\in J\} \\
    &\qquad\quad\sqcup \{v\leftarrow f_j, v'\leftarrow f_j, v_j\leftarrow f_j:j\in J\}\sqcup \{ e\leftarrow f_j, e_j\leftarrow f_j, e'_j\leftarrow f_j: j\in J\}
\end{align*}
Figure \ref{fig:cYfor32} illustrates this scwol for $q=3$.
\begin{figure}
\centering
\begin{tikzpicture}[scale=4,x=30,y=30,z=10]
    \makeatletter
    \define@key{cylindricalkeys}{angle}{\def\myangle{#1}}
    \define@key{cylindricalkeys}{radius}{\def\myradius{#1}}
    \define@key{cylindricalkeys}{z}{\def\myz{#1}}
    \tikzdeclarecoordinatesystem{cylindrical}
    {
       \setkeys{cylindricalkeys}{#1}
       \pgfpointadd{\pgfpointxyz{0}{0}{\myz}}{\pgfpointpolarxy{\myangle}{\myradius}}
    }

    \tikzstyle{every edge}=[draw,->]
    \node (vv) at (0,0,0) {$v'$};
    \node (v) at (0,0,2) {$v$};
    \node (e) at (0,0,1) {$e$} edge (v) edge (vv);

    \foreach \a in {4,3,2,1,0}
    {
        \ifnum \a > 1
            \tikzstyle{every node}=[color=gray];
            \tikzstyle{every edge}=[draw,->,color=lightgray];
        \fi

        \node (v1) at (cylindrical cs:angle=\a*72+20,radius=1,z=1) {$v_\a$};
        \node (e1) at (cylindrical cs:angle=\a*72+20,radius=.5,z=1.5) {$e_\a$} edge (v) edge (v1);
        \node (ee1) at (cylindrical cs:angle=\a*72+20,radius=.5,z=.5) {$e'_\a$} edge (vv) edge (v1);

        \ifnum \a > 0 
        \node (f1) at (cylindrical cs:angle=\a*72+20,radius=.333,z=1) {$f_\a$} edge (v) edge (vv) edge (e) edge (e1) edge (ee1) edge (v1);
        \else
        \node (f1) at (cylindrical cs:angle=\a*72+20,radius=.333,z=1) {$f_\a$} edge (v) edge (vv) edge (e) edge (e1) edge (v1);
        \draw[->] (ee1) -- (f1);
        \fi

    }

\end{tikzpicture} \caption{The scwol $\cY$ for $q=3$}\label{fig:cYfor32}
\end{figure}

We construct a complex of groups $G(\cY)$ over $\cY$ by setting the vertex groups to be
\[
G_v = S,\qquad G_{v'}=S',\qquad G_{v_j}= S_{\lambda(j)} \times S'_{\lambda'(j)},\qquad G_{e_j}=S_{\lambda(j)},\qquad G_{e'_j}=S'_{\lambda'(j)}.
\]
All other vertex groups are trivial. The inclusions are the obvious ones. All twist elements are trivial, so $G(\cY)$ is a \emph{simple} complex of groups.

We endow $|\cY|$ with a locally Euclidean metric as follows: Let $\Delta$ be the geometric realisation of one triangle in the affine Coxeter complex of type $\tilde C_2$. For each $j\in J$, we map each subcomplex spanned by $\{v,v',v_j,e,e_j,e'_j,f_j\}$ onto the barycentric subdivision of $\Delta$ in the obvious way and pull back the metric. We obtain a locally Euclidean metric on $|\cY|$ with angle $\pi/4$ at the vertices $v,v'$ and with angle $\pi/2$ at the vertex $v_j$.

\begin{Lemma}\label{la:c2_local_development}
	The local development $\cY(\tilde v_j)$ is isomorphic to the cone over the barycentric subdivision of a complete bipartite graph of order $(q,q)$. The geometric link $\lk(\tilde v_j,\st(\tilde v_j))$ is hence isometric to the barycentric subdivision of a generalised $2$-gon, in particular a connected, CAT(1) polyhedral complex of diameter $\pi$.
\end{Lemma}

\begin{Proof}
	We abbreviate $H=S_{\lambda(j)}$, $F=S'_{\lambda'(j)}$ and $G=H\times F$. We have $\cY(\tilde v) = \{\tilde v\} * \Lk_{\tilde v_j}(\cY)$ by the definition of the local development, where
	\begin{align*}
		V(\Lk_{\tilde v_j}(\cY)) &= \{ (gH, e_j) : g\in G \} \sqcup \{ (gF,e'_j) : g\in G\} \sqcup \{ (g,f_j) : g \in G \} \cong F \sqcup H \sqcup G\\
		E(\Lk_{\tilde v_j}(\cY)) &= \{ (g, e_j\leftarrow f_j) : g\in G\} \sqcup  \{ (g, e'_j\leftarrow f_j) : g\in G\},
	\end{align*}
    which is easily seen to be the barycentric subdivision of a complete bipartite graph on the vertex set $F\sqcup H$. The rest of the argument is analogous to the proof of Proposition \ref{prop:key}.
\end{Proof}

\begin{Theorem}\label{th:c2_result2}
    The complex $G(\cY)$ is developable. The fundamental group $\Gamma=\pi_1(G(\cY))$ admits the presentation
    \[
    \Gamma = \langle S,S' \,|\, \text{ all relations in the groups } S,S',\; [S_{\lambda(j)},S'_{\lambda'(j)}]=1\quad \forall j\in J \rangle.
    \]
    The universal cover $\cX$ of this complex of groups is a building of type $\tilde C_2$, and $\Gamma$ is a uniform lattice acting regularly on one type of panels of $\cX$.
\end{Theorem}

\begin{Proof}
    This proof is analogous to the proofs of Proposition \ref{prop:c2_developable} and Theorem \ref{th:c2_building}, where we use Lemma \ref{la:c2_local_development} for the local developments at the vertices $v_j$.
\end{Proof}

\noindent In the next two small sections, we will make this explicit in the cases where we have simple presentations.

\subsubsection{The prime case}

If $q=p$ is prime and if we start with two slanted symplectic quadrangles $\cI$, $\cI'$ of order $(p-1,p+1)$ with corresponding Heisenberg groups $S$, $S'$ with respective generator pairs $x$, $y$ and $x'$, $y'$, the sets of line representatives $L$ and $L'$ can be parametrised by $\bP\F_p^2\sqcup\{0\}$. The stabilisers then have the form 
\[
S_{[a:b]}=\langle x^ay^bz^{-\tfrac{1}{2}ab}\rangle\qquad\text{and}\qquad S_0=\langle z\rangle,
\]
where $z=[x,y]$ and similarly for $S'$. Let $J=\{0,1,\ldots, p+1\}$.

\begin{Theorem}\label{th:c2_prime_result2}
	Let $\lambda, \lambda':J\rightarrow \bP\F_p^2\sqcup \{0\}$ be two bijections. The group
    \[
    \Gamma = \Biggl\langle x,y,x',y' \,\Bigg|\, \begin{array}{c} z= xyx^{-1}y^{-1}, z'=x'y'x'^{-1}y'^{-1}, \\
            x^p=y^p=x'^p=y'^p=z^p=z'^p=1, \\ xz=zx, yz=zy, x'z'=z'x', y'z' = z'y', \\\space
            [S_{\lambda(j)}, S'_{\lambda'(j)}] = 1\qquad \forall j\in J \end{array}\Biggr\rangle
    \]
    is a uniform lattice in a building of type $\tilde C_2$.
\end{Theorem}

\begin{Proof}
	This is a direct application of Theorem \ref{th:c2_result2}.
\end{Proof}

\subsubsection{The abelian case}

If $q$ is even, then the associated Singer groups $S$ and $S'$ are isomorphic to $\F_q^3$. In this case, we obtain a very simple presentation:

Let $J=\{0,1,\ldots,q+1\}$ and fix two bijections $\lambda,\lambda': J \rightarrow \bP\F_q^2\sqcup \{0\}$. We write $S_{[a:b]}=\F_q(a,b,0)^T$ and $S_{0}= \F_q(0,0,1)^T$ for the line stabilisers in $S$ and add a prime $'$ for the respective groups in $S'$.

\begin{Theorem}
	The group
	\[
	\Gamma = \bigl(S * S' \bigr) / \langle [ S_{\lambda(j)}, S'_{\lambda'(j)}] : j\in J\rangle
	\]
	is a uniform lattice in an exotic building of type $\tilde C_2$.
\end{Theorem}

\begin{Proof}
    Again, this is a direct application of Theorem \ref{th:c2_result2}.
\end{Proof}

\begin{Remark}
    Note that this is a finitely presented group generated by involutions such that all other relations make elements commute. It would be very interesting if one of these groups were a right-angled Coxeter group, but this seems very unlikely.
\end{Remark}

\subsection{Property (T)}

Even though most of these buildings are necessarily exotic, the lattices have automatically property (T).

\begin{Proposition}[\.Zuk]
	All cocompact lattices in buildings of type $\tilde C_2$ have property (T) if the order of the generalised quadrangles which appear as vertex links is not $(2,2)$. In particular, all lattices we construct here have property (T).
\end{Proposition}

\begin{Proof}
	See \cite{Zuk:TGP:96}. A more detailed exposition can be found in \cite{BHV:KPT:08}.
\end{Proof}

\section{Group homology}\label{sec:group_homology}

We will calculate group homology of the lattices we have constructed using the action on the building. Remember that a group action on a polyhedral complex is said to be \emph{without inversions} if every setwise stabiliser of a polyhedral cell stabilises the cell pointwise. The actions of the lattices we have constructed are obviously without inversions, since the lattices act in a type-preserving fashion.

\begin{Theorem}[VII.7 and VII.8 in \cite{Bro:CoG:82}]\label{th:spec_seq}
    If a group $G$ acts without inversions on a contractible polyhedral complex $X$, and if $\Sigma_i$ is a system of representatives for the $i$-cells of $X$, there is a spectral sequence
    \[
    E^1_{i,j}= \bigoplus_{\sigma\in\Sigma_i} H_j(G_{\sigma}) \Rightarrow H_{i+j}(G).
    \]
    In particular, we have $E^1_{i,0}\cong \bigoplus_{\sigma\in\Sigma_i} \Z$ and the differential $d^1_{i,0}$ is induced by the differential $\partial_{i}$ on cellular chains. This results in
    \[
    E^2_{i,0}\cong H_i(G\backslash X).
    \]
\end{Theorem}

\begin{Proof}
    The construction of the spectral sequence can be found in \cite[VII.7]{Bro:CoG:82}. The structure of the differential $d^1$ on the bottom row can either be seen directly from the construction of the spectral sequence as the spectral sequence associated to group homology with coefficients in cellular chains or by specialising \cite[VII.8]{Bro:CoG:82}.
\end{Proof}

\begin{Corollary}\label{cor:rational_homology}
    If all stabilisers $G_\sigma$ are finite for every cell $\sigma\in X$, we have \[H_*(G;\Q)\cong H_*(G\backslash X;\Q).\]
\end{Corollary}

\begin{Proof}
    This follows easily from Theorem \ref{th:spec_seq} and the fact that $H_j(G_\sigma;\Q)=0$ for all $j>0$ and for all simplices $\sigma$ by the transfer map, see \cite[Corollary III.10.2]{Bro:CoG:82}.
\end{Proof}

Since cyclic groups are abundant in the construction of our lattices, we need the following simple result.

\begin{Proposition}[II.3.1 in \cite{Bro:CoG:82}]\label{prop:hom_of_cycl_groups}
    Let $G$ be a cyclic group of order $n$. Then
    \[
        H_j(G) \cong\begin{cases}
            \Z & j=0 \\
            \Z/n & j\text{ odd}\\
            0 & \text{ else.}
        \end{cases}
    \]
\end{Proposition}

\subsection{Group homology of lattices in buildings of type \texorpdfstring{$\tilde A_2$}{\textasciitilde A2}}

We will calculate the full group homology of the lattices generated by cyclic Singer groups which we have constructed in section \ref{subsec:cyclic_lattices}. The abelianisation can be calculated directly.

\begin{Proposition}
    Let $\Gamma_1$ be a lattice constructed as in Theorem \ref{th:cyclic_lattices} using a classical projective plane of order $q$ and three ordered classical difference sets $\delta_1$, $\delta_2$ and $\delta_3$. We have
    \[
    H_1(\Gamma_1) = \Bigl\langle s_1,s_2,s_3 \,\Big|\, \begin{array}{c} (q^2+q+1)s_1 = (q^2+q+1)s_2 = (q^2+q+1)s_3 = 0 \\ \delta_1(j)s_1+\delta_2(j)s_2+\delta_3(j)s_3=0\quad\forall j\in J \end{array}\Bigr\rangle.
    \]
    In particular, $\Gamma_1$ is perfect if and only if the $(q\times 3)$-matrix $\cD= (\delta_j(i))_{i,j}$ over $\Z/(q^2+q+1)$ has full rank in the sense of \cite[Chapter 4]{Bro:MCR:93}.
\end{Proposition}

\begin{Proof}
    Since $H_1(\Gamma_1)$ is isomorphic to the abelianisation of $\Gamma_1$, we just have to abelianise the presentation from Theorem \ref{th:cyclic_lattices}. The resulting group is the kernel of the linear map $(\Z / (q^2+q+1))^3 \rightarrow (\Z / (q^2+q+1))^q$ described by $\cD$. In particular, by \cite[Theorem 5.3]{Bro:MCR:93}, the kernel is trivial if and only if $\cD$ has full rank.
\end{Proof}

\begin{Theorem}
    Again, let $\Gamma_1$ be a lattice constructed as in Theorem \ref{th:cyclic_lattices}. Let $q$ be the order of the associated projective plane. Then
    \[
        H_j(\Gamma_1) \cong\begin{cases}
            \Z & j=0\\
            \ker(\cD) & j=1 \\
            \Z^q & j=2 \\
            (\Z/(q^2+q+1))^3 & j\geq 3 \text{ odd} \\
            0 & \text{ else.}
        \end{cases}
    \]
    In addition
    \[
    H_2(\Gamma_2;\Q)=\Q^q,\qquad H_j(\Gamma_2;\Q)=0 \text{ for }j\not\in\{ 0,2\}
    \]
    for any lattice $\Gamma_2$ constructed as in Theorems \ref{th:a2_is_building} or \ref{th:cyclic_lattices}.
\end{Theorem}

\begin{Proof}
    Consider the spectral sequence $E^1_{i,j}$ from Theorem \ref{th:spec_seq} for the $\Gamma_1$-action on the associated building $X$. The quotient space $\Gamma_1\backslash X$ is the geometric realisation of the scwol $\cY$ from the construction in section \ref{subsec:general_a2_construction}, which has hence the homotopy type of a bouquet of $2$-spheres. The structure of the spectral sequence is particularly simple since most stabiliser subgroups are trivial. We have
    \[
    E^1_{0,j} \cong \bigoplus_{k=1}^3 H_j( \Z/(q^2+q+1)), \qquad E^1_{1,0} \cong \Z^3, \qquad E^1_{2,0} \cong \Z^{q+1}.
    \]
    All other groups on the first page are trivial. The differential on the bottom row is induced by the cellular differential, so we have
    \[
    E^2_{0,0}\cong \Z,\qquad E^2_{1,0}=0, \qquad E^2_{2,0}\cong \Z^q.
    \]
    All other groups remain unchanged. Since the only remaining nontrivial differential is
    \[
    d^2_{2,0}: E^2_{2,0}\cong \Z^q \rightarrow (\Z/(q^2+q+1))^3 = E^2_{0,1},
    \]
    whose kernel is always isomorphic to $\Z^q$, we obtain the full description of $H_*(\Gamma_1)$ by inspecting the third page $E^3_{i,j}$ and using Proposition \ref{prop:hom_of_cycl_groups} for the homology of cyclic groups. Rational homology can easily be calculated using Corollary \ref{cor:rational_homology}.
\end{Proof}

\subsection{Group homology of lattices in buildings of type \texorpdfstring{$\tilde C_2$}{\textasciitilde C2}}

We will first calculate rational group homology for the lattices of type $\tilde C_2$ we have constructed in section \ref{sec:c2}.

\begin{Theorem}\label{th:c2_rational_homology}
    For a lattice $\Gamma$ constructed as in sections \ref{subsec:c2_two_panels} or \ref{subsec:c2_one_panel}, we have $H_j(\Gamma;\Q)=0$ for $j\neq 0$.
\end{Theorem}

\begin{Proof}
    The geometric realisations of the quotient scwols are contractible in both cases. We can hence apply Corollary \ref{cor:rational_homology} to obtain the result.
\end{Proof}

\noindent For the lattices acting regularly on two types of panels as in section \ref{subsec:c2_two_panels}, the situation with integral coefficients is rather complicated and probably depends on the identification bijections $\lambda$ and $\lambda'$ given in the construction.

Hence we will consider lattices acting regularly on one type of panel as in section \ref{subsec:c2_one_panel}. Assume that the associated slanted symplectic quadrangle has order $(q-1,q+1)$, denote the associated Singer group by $S$.

\begin{Theorem}
    For the lattice $\Gamma$, we have
    \[
    H_1(\Gamma)\cong (\Z/q)^6\qquad\text{and}\qquad H_2(\Gamma) \cong H_2(S) \oplus H_2(S).
    \]
\end{Theorem}

\begin{Proof}
    The structure of the abelianisation can easily be seen using the presentation from Theorem \ref{th:c2_result2}. For $H_2$, we inspect the corresponding spectral sequence from Theorem \ref{th:spec_seq}. As a set of representatives for the action, we use the subcomplex induced by the set of all chambers containing a common panel of the type the lattice acts regularly on. Then we have
    \begin{align*}
        E^1_{0,2} &\cong H_2(S)^2,& E^1_{1,2} &=0,\\
        E^1_{0,1} &\cong (\Z/q)^6 \oplus (\Z/q)^{2q+4},& E^1_{1,1} &\cong (\Z/q)^{2q+4},\\
        E^1_{0,0} &\cong \Z^{q+4},& E^1_{1,0}&\cong \Z^{2q+5},& E^1_{2,0}&\cong \Z^{q+2}.
    \end{align*}
    In this calculation, we use Proposition \ref{prop:hom_of_cycl_groups} as well as $H_1(S)=(\Z/q)^3$, which is easily verified. Since the quotient space $\Gamma\backslash X$ is contractible, we obtain
    \[
    E^2_{0,0} \cong \Z,\qquad E^2_{1,0}=0,\qquad E^2_{2,0}=0.
    \]
    Since we already know that $H_1(\Gamma)\cong (\Z/q)^6$, the map $d^1_{1,1}:E^1_{1,1}\rightarrow E^1_{0,1}$ must be injective and we obtain $E^2_{1,1}=0$ which proves the theorem.
\end{Proof}

\nocite{BL:TL:01}
\bibliographystyle{alpha}
\bibliography{../../biblio}

\begin{thebibliography}{CMSZ93b}

\bibitem[BdlHV08]{BHV:KPT:08}
B.~Bekka, P.~de~la Harpe, and A.~Valette.
\newblock {\em Kazhdan's property ({T})}, volume~11 of {\em New Mathematical
  Monographs}.
\newblock Cambridge University Press, Cambridge, 2008.

\bibitem[BH99]{BH:NPC:99}
M.~R. Bridson and A.~Haefliger.
\newblock {\em Metric spaces of non-positive curvature}, volume 319 of {\em
  Grundlehren der Mathematischen Wissenschaften [Fundamental Principles of
  Mathematical Sciences]}.
\newblock Springer-Verlag, Berlin, 1999.

\bibitem[BL01]{BL:TL:01}
H.~Bass and A.~Lubotzky.
\newblock {\em Tree lattices}, volume 176 of {\em Progress in Mathematics}.
\newblock Birkh\"auser Boston Inc., Boston, MA, 2001.

\bibitem[Bro82]{Bro:CoG:82}
K.~S. Brown.
\newblock {\em Cohomology of groups}, volume~87 of {\em Graduate Texts in
  Mathematics}.
\newblock Springer-Verlag, New York, 1982.

\bibitem[Bro89]{Bro:Bdg:89}
K.~S. Brown.
\newblock {\em Buildings}.
\newblock Springer-Verlag, New York, 1989.

\bibitem[Bro93]{Bro:MCR:93}
W.~C. Brown.
\newblock {\em Matrices over commutative rings}, volume 169 of {\em Monographs
  and Textbooks in Pure and Applied Mathematics}.
\newblock Marcel Dekker Inc., New York, 1993.

\bibitem[CL01]{CL:MC:01}
R.~Charney and A.~Lytchak.
\newblock Metric characterizations of spherical and {E}uclidean buildings.
\newblock {\em Geom. Topol.}, 5:521--550 (electronic), 2001.

\bibitem[CMS94]{CMS:TA2:94}
D.~I. Cartwright, W.~M{\l}otkowski, and T.~Steger.
\newblock Property ({T}) and {$\tilde{A\sb 2}$} groups.
\newblock {\em Ann. Inst. Fourier (Grenoble)}, 44(1):213--248, 1994.

\bibitem[CMSZ93a]{CMSZ1}
D.~I. Cartwright, A.~Mantero, T.~Steger, and A.~Zappa.
\newblock Groups acting simply transitively on the vertices of a building of
  type {$\tilde A\sb 2$}. {I}.
\newblock {\em Geom. Dedicata}, 47(2):143--166, 1993.

\bibitem[CMSZ93b]{CMSZ2}
D.~I. Cartwright, A.~Mantero, T.~Steger, and A.~Zappa.
\newblock Groups acting simply transitively on the vertices of a building of
  type {$\tilde A\sb 2$}. {II}. {T}he cases {$q=2$} and {$q=3$}.
\newblock {\em Geom. Dedicata}, 47(2):167--223, 1993.

\bibitem[Dem68]{De:FG:68}
P.~Dembowski.
\newblock {\em Finite geometries}.
\newblock Ergebnisse der Mathematik und ihrer Grenzgebiete, Band 44.
  Springer-Verlag, Berlin, 1968.

\bibitem[FHT09]{FHT:PAG:09}
B.~Farb, C.~Hruska, and A.~Thomas.
\newblock Problems on automorphism groups of nonpositively curved polyhedral
  complexes and their lattices.
\newblock To appear in Geometry, Topology and Rigidity (B. Farb and D. Fisher,
  eds.), 2009.

\bibitem[GAP]{GAP4}
The GAP~Group.
\newblock {\em {GAP -- Groups, Algorithms, and Programming, Version 4.4.12}}.
\newblock (\url{http://www.gap-system.org}).

\bibitem[Gil07]{Gil:TPP:07}
N.~Gill.
\newblock Transitive projective planes.
\newblock {\em Adv. Geom.}, 7(4):475--528, 2007.

\bibitem[GJS94]{GJS:SSQ:94}
T.~Grundh{\"o}fer, M.~Joswig, and M.~Stroppel.
\newblock Slanted symplectic quadrangles.
\newblock {\em Geom. Dedicata}, 49(2):143--154, 1994.

\bibitem[HvM89]{HvM:PHP:89}
G.~Hanssens and H.~van Maldeghem.
\newblock On projective {H}jelmslev planes of level {$n$}.
\newblock {\em Glasgow Math. J.}, 31(3):257--261, 1989.

\bibitem[Kan86]{Kan:GPS:86}
W.~M. Kantor.
\newblock Generalized polygons, {SCAB}s and {GAB}s.
\newblock In {\em Buildings and the geometry of diagrams ({C}omo, 1984)},
  volume 1181 of {\em Lecture Notes in Math.}, pages 79--158. Springer, Berlin,
  1986.

\bibitem[KLT87]{KLT:CTL:87}
W.~M. Kantor, R.~A. Liebler, and J.~Tits.
\newblock On discrete chamber-transitive automorphism groups of affine
  buildings.
\newblock {\em Bull. Amer. Math. Soc. (N.S.)}, 16(1):129--133, 1987.

\bibitem[KMW84]{KMW:A2l:84}
P.~K{\"o}hler, T.~Meixner, and M.~Wester.
\newblock The affine building of type {$\tilde A\sb{2}$} over a local field of
  characteristic two.
\newblock {\em Arch. Math. (Basel)}, 42(5):400--407, 1984.

\bibitem[LaJ]{LaJolla}
D. Gordon.
\newblock {\em La {J}olla difference set repository}.
\newblock (\url{http://www.ccrwest.org/diffsets/diff\_sets/}).

\bibitem[Mar91]{Mar:DSS:91}
G.~A. Margulis.
\newblock {\em Discrete subgroups of semisimple {L}ie groups}, volume~17 of
  {\em Ergebnisse der Mathematik und ihrer Grenzgebiete (3) [Results in
  Mathematics and Related Areas (3)]}.
\newblock Springer-Verlag, Berlin, 1991.

\bibitem[Ron84]{Ron:TG:84}
M.~A. Ronan.
\newblock Triangle geometries.
\newblock {\em J. Combin. Theory Ser. A}, 37(3):294--319, 1984.

\bibitem[Ron86]{Ron:CBR:86}
M.~A. Ronan.
\newblock A construction of buildings with no rank {$3$} residues of spherical
  type.
\newblock In {\em Buildings and the geometry of diagrams ({C}omo, 1984)},
  volume 1181 of {\em Lecture Notes in Math.}, pages 242--248. Springer,
  Berlin, 1986.

\bibitem[Ron89]{Ron:LoB:89}
M.~A. Ronan.
\newblock {\em Lectures on buildings}, volume~7 of {\em Perspectives in
  Mathematics}.
\newblock Academic Press Inc., Boston, MA, 1989.

\bibitem[SAG]{SAGE}
SAGE~Mathematical Software.
\newblock {\em {Version 4.1}}.
\newblock (\url{http://www.sagemath.org}).

\bibitem[Sin38]{Si:FPG:38}
J.~Singer.
\newblock A theorem in finite projective geometry and some applications to
  number theory.
\newblock {\em Trans. Amer. Math. Soc.}, 43(3):377--385, 1938.

\bibitem[STdW09]{STW:SQ:09}
E.~E. Shult, K.~Thas, and S.~de~Winter.
\newblock {\em Singer Quadrangles}.
\newblock Oberwolfach Preprint OWP 2009-07. MFO Oberwolfach, 2009.

\bibitem[Str03]{Str:PSQ:03}
M.~Stroppel.
\newblock Polarities of symplectic quadrangles.
\newblock {\em Bull. Belg. Math. Soc. Simon Stevin}, 10(3):437--449, 2003.

\bibitem[Tit74]{Tit:BsT:74}
J.~Tits.
\newblock {\em Buildings of spherical type and finite {BN}-pairs}.
\newblock Springer-Verlag, Berlin, 1974.
\newblock Lecture Notes in Mathematics, Vol. 386.

\bibitem[vM90]{HvM:NTB:90}
H.~van Maldeghem.
\newblock Automorphisms of nonclassical triangle buildings.
\newblock {\em Bull. Soc. Math. Belg. S\'er. B}, 42(2):201--237, 1990.

\bibitem[Wei09]{Wei:SAB:09}
R.~M. Weiss.
\newblock {\em The structure of affine buildings}, volume 168 of {\em Annals of
  Mathematics Studies}.
\newblock Princeton University Press, Princeton, NJ, 2009.

\bibitem[{\.Z}uk96]{Zuk:TGP:96}
A.~{\.Z}uk.
\newblock La propri\'et\'e ({T}) de {K}azhdan pour les groupes agissant sur les
  poly\`edres.
\newblock {\em C. R. Acad. Sci. Paris S\'er. I Math.}, 323(5):453--458, 1996.

\end{thebibliography}
\end{document}